\documentclass[12pt,a4paper]{article}
\usepackage[left=2cm,right=2cm,top=3cm,bottom=3cm]{geometry}
\usepackage{graphicx} 
\usepackage{dutchcal}
\graphicspath{ {C:\Latex\TWS990} }
\usepackage{booktabs} 
\usepackage{amssymb}
\usepackage{amsmath}
\numberwithin{equation}{section}
\usepackage{amsthm}
\usepackage{textcomp}
\usepackage{color}
\usepackage{xcolor}
\usepackage{mathrsfs}
\newtheorem{theorem}{Theorem}[section]
\newtheorem{corollary}{Corollary}[section]

\newtheorem{example}{Example}[section]
\newtheorem{lemma}{Lemma}[section]
\newtheorem{proposition}{Proposition}[section]
\newtheorem{remark}{Remark}[section]
\usepackage{url}
\usepackage[utf8]{inputenc}
\usepackage[english]{babel}

\newcommand{\mc}[1]{\mathcal{#1}}
\newcommand{\mbb}[1]{\mathbb{#1}}
\newcommand{\mb}[1]{\boldsymbol{#1}}
\newcommand{\mr}[1]{\mathring{#1}}
\newcommand{\la}{\lambda}
\newcommand{\cl}[2]{\int_{#1}^{#2}}
\newcommand{\ov}[1]{\overline{#1}}

\newcommand{\mf}[1]{\mathfrak{#1}}
\newcommand{\e}{\varepsilon}
\DeclareMathOperator*{\esssup}{ess\,sup}
\newcommand{\sem}[1]{({#1}(t))_{t\geq 0}}
\newcommand{\p}{\partial}
\title{Fragmentation--coagulation processes with advection or diffusion in space}
\author{Jacek Banasiak\footnote{Corresponding author: jacek.banasiak@up.ac.za} \footnote{The research of J.B. and N.M.  was supported by the National Research Foundation of South Africa, grant 87720. J.B. was partially supported by the IMPRESS-U programme in the framework of the project
\#2023/05/Y/ST6/00263} \footnote{The paper was completed during J. Banasiak's fellowship at the Stellenbosch Institute for Advanced Studies.}
 \\
\small{University of Pretoria, Łódź University of Technology}\\\& \small{Stellenbosch Institute for Advanced Studies}\\
Nduduzo Majozi\\
\small{University of Pretoria}}
\date{}

\begin{document}
\maketitle
\begin{abstract}
\noindent
In this paper, we consider a continuous fragmentation--coagulation model in which the reacting particles can be transported in physical space through either advection or diffusion.  We prove new results on the generation of $C_0$-semigroups with parameter and use them to show that the Abstract Cauchy Problem associated with a more general version of the advection/diffusion--fragmentation problem generates a positive $C_0$-semigroup in spaces $L_1(\mbb R_+, X_x, (1+m^r)dm),$ where $m$ is the particle mass, $X_x$ is the space of either integrable or continuous functions with respect to the spatial variable, and the weight exponent $r$ is sufficiently large. These results enable us to prove the classical solvability of a wide range of advection/diffusion--fragmentation--coagulation equations with unbounded coagulation kernels.  \end{abstract}

		\section{Introduction}
        \subsection{Fragmentation and coagulation processes}
	Fragmentation and coagulation are processes in which animate or inanimate objects split into smaller parts or, conversely, smaller objects join together to form a larger one. Such processes occur in various fields, from polymer science, where we observe polimerisation and depolimerisation of molecules, to population dynamics, where animals or humans form larger groups or split into smaller ones in response to environmental changes. They are also important in industrial practice, e.g., in water treatment, suspended particles are removed from the liquid by adding coagulants, which neutralise charges and allow particles to clump together and sink to the bottom, and fragmentation is used to produce powders of uniformly-sized grains. A comprehensive overview of the field can be found in \cite[Section 2]{BLL1}.  Fragmentation and coagulation processes are often accompanied by internal dynamics, such as the growth or decay of clusters due to birth or death when clusters are formed by animate matter, or, in chemical applications, by material deposition and dissolution. Such processes have received substantial coverage in the literature, see, e.g., \cite{Bana12a, BaAr, BaLa09, BanLam20, BaPoSh22, CEH91} and \cite[Section 5.2]{BLL1}, or \cite{BaPiRu, BerGab20, DoGa10, MKJB, PerTr}, where the emphasis is on the long-term behaviour of the dynamics. A discrete version of such models has been analysed comprehensively in \cite{Banasiak2019, Banasiak2018}. 
    
    The above results assume that the processes occur in a homogeneous medium, so there is no explicit dependence on the spatial variable. If, however, the dynamics in physical space is relevant, the coefficients in the models must be made space-dependent, and the equations complemented by terms describing spatial transport.  
    
    \subsection{The model}
    
    To introduce our model, we assume that we have an ensemble of clusters (also called particles) that, at a given time $t$, can be characterised by their position in space $x\in \mbb R^N$ and another scalar continuous attribute $m\in \mbb R_+:=[0,\infty)$ (in natural sciences, typically $m$ is the mass or size of the cluster). The state of the ensemble at time $t$ is described by the density $u(t,x,m)$ of clusters at location $x$ and mass $m$. We consider clusters to be material points of mass $m$ that can split into smaller clusters or combine with another cluster to form a bigger one (we neglect possible coagulation of more than two clusters). In general, both fragmentation and coagulation may be nonlocal in space, meaning that the daughter clusters can appear at different positions from the parent particle, and the coagulating clusters could initially occupy different positions than the combined one; see, e.g., \cite{LaWr00, RudWiecz1}. Considering such a scenario would make the model similar to the Enskog equation, where the interactions between particles occur on their surfaces and thus the volume occupied by them is excluded, or the Povzner equation, where the interactions are allowed to occur also within it, \cite{ bellach, ens, Pov}; the papers \cite{LaWr00, RudWiecz1} are concerned with this approach. However, following the subject literature discussed below, we simplify the model (but not necessarily the mathematics) by considering what in the kinetic theory is referred to as the Boltzmann limit and assume that the action occurs locally in space, that is, the daughter clusters remain at the same point as the parent cluster, and the coagulation occurs between clusters occupying the same point in space. 

    In this paper, we consider the motion of particles in the surrounding medium, either by diffusion or by advection, the latter arising from the flow of the medium. When it does not lead to confusion, we will use the term transport to refer to both processes.

    Thus, the governing equation for the density $u$ takes the form   
    \begin{subequations}\label{GeneralEqn}
    \begin{equation}\label{GeneralEqna}
\partial_tu(t,x,m)=\mathcal{T}_0u(t,x,m)+\mathcal{F}u(t,x,m)+\mathcal{C}u(t,x,m), \hspace{.5cm} t>0,~(x,m)\in\Omega\times \mathbb{R}_+,
	\end{equation} 
	and is supplemented with the initial condition \begin{equation}\label{ic1}
		u(0,x,m)=\mr{u}(x,m), \quad (x,m)\in\Omega\times\mathbb{R}_+,
	\end{equation}
    \end{subequations}
	where $\Omega\subseteq \mbb R^N,$ $N\geq 1,$ is an open set, $\mbb R_+=[0,\infty),$ $\mc T_0$ is a family of general transport expressions depending on $m\in \mbb R_+$ as a parameter, $(\mc T_{0,m})_{m\in \mbb R_+},$ acting as \begin{equation}
(x,m) \mapsto [\mc T_{0,m}u(\cdot,m)](x),\label{T0}
\end{equation}
 and the fragmentation and coagulation processes are described by 
	\begin{subequations}\label{ops}
                    \begin{equation}\begin{split}
		\mathcal{F}u(x,m)= \mc Au(x,m)+\mc Bu(x,m)=-a(x,m)u(x,m)\!+\!\int_{m}^{\infty}\!\!b(x,m,s)a(x,s)u(x,s)ds,\end{split}\label{ops2}
        \end{equation}
        and 
        \begin{equation}
        \mc Cu(x,m):= [\ov{\mc C}(u,u)](x,m),
        \label{mcCC}
        \end{equation}
        \end{subequations}
        where
        $$
    [\ov{\mathcal{C}}(u,v)](x,m)=\dfrac{1}{2}\int_{0}^{m}k(x, m-s,s)u(x,m-s)v(x,s)ds-u(x,m)\int_{0}^{\infty}k(x, m,s)v(x,s)ds.
        $$
	In the equations above, $\mc F=\mc A+\mc B$ is referred to as the fragmentation operator, split into the so-called loss operator $\mc A$ and the gain operator $\mc B$, with the fragmentation rate $a $  and
 the fragmentation kernel $ b $ that describes the distribution of masses  $m$ of the daughter particles at position $x,$ spawned by splitting of a particle of mass $s$; $\mc C$ denotes the coagulation operator with a positive and symmetric in $s$ and $m$ function $k$ giving the rate of coagulation, see \cite{BLL1}. 

Most work on \eqref{GeneralEqn} focuses on matter diffusing through space, see, e.g., the derivation from an individual-based model in \cite{RudWiecz1} (also \cite[Section 7.10.1]{Oku}). The weak solvability of the problem with diffusion and discrete coagulation and/or fragmentation was considered in a series of papers \cite{BeWr97, LaWr00, LaMi02c, LaWr98a, LaWr98b, Wrzo97, Wrzo02}, and with continuous ones in \cite{Buro85, LaMi02a}. Considerable attention was devoted to long-term dynamics of the problem \cite{DeFe09, DeFe14, MiRR03}. Classical solutions for diffusion with continuous fragmentation and bounded coagulation were obtained using the theory of analytic semigroups in \cite{Aman00, AmWa05}, and in \cite{shindin2026, Walk05} these results were extended to an unbounded coagulation operator;  in the former, the diffusion and the coagulation processes should be suitably dominated by the fragmentation, and in the latter, the fragmentation and coagulation are controlled by the diffusion term.  In particular, \cite{shindin2026} provides the existence of global in time solutions for a class of such problems.   
 A survey of results on diffusion--fragmentation--coagulation, with more comprehensive references, can be found in \cite[Section 11.2]{BLL2}. The case in which the particle transport is described by a first-order partial differential operator was considered in, e.g., \cite{Buro83, ChDu97, Dubo90, Dubo01}, where the authors used successive approximation techniques in spaces with an exponential weight for the integral version of the problem, written in characteristic coordinates. We also mention \cite{Clov}, which introduces a framework based on techniques from \cite{BanTTSP2007} but contains several gaps. Recently, such processes have been considered from another viewpoint in \cite{ChrisNietVel, CrisVel}.

 In this paper,  we provide a comprehensive approach to the problem based on the theory of $C_0$-semigroups in the spirit of \cite{Bana12b, BanLam20}, that is, we prove the existence and moment regularising property of a $C_0$-semigroup generated by the linear part of the equation. Here, we introduce a novel approach, constructing an equation that dominates the transport--fragmentation part of the original equation, but in which the new fragmentation operator is independent of the spatial variable $x$. If this new fragmentation kernel is uniformly integrable with respect to the parent's size (e.g., homogeneous fragmentation kernels have this property), then there is $r_1$ such that for all $r\geq r_1$ the original transport--fragmentation operator generates a $C_0$-semigroup in $\mc X^1_r:= L_1(\mbb R_+\times \Omega, dm_rdx), dm_r := (1+m^r)dm$.  We note that similar results for the pure fragmentation operator have been recently obtained in \cite{kerr2025c, kerr2025}, but the authors constructed there general (exponential) weights, whereas here we remain within the polynomial ones.  
 
 The full equation is treated as a semilinear perturbation of the transport--fragmentation part, and the moment regularisation property of the latter allows us to prove the classical solvability of problems with unbounded coagulation, as long as it grows slower than the fragmentation rate as particles become large. We emphasise that, in contrast to \cite{Walk05}, where the coagulation is controlled by the growth rate of the diffusion coefficient, this paper uses a regularising effect of the loss term.   
 
 We note, however, that while the results for the linear transport--fragmentation part are valid in $\mc X^1_r$  (physically natural, as for $r=1$ the norm gives the total mass and number of particles of the system), passing to the nonlinear model encounters a serious difficulty because we have chosen to work with a model that is local with respect to the space variable, see \cite[Section 11.2]{BLL2}. Thus, staying within the framework of $L_1$ spaces, we encounter a product of two $L_1$ functions in the coagulation term, which, in general, is not well-defined in  $L_1$. The references above provide various methods for addressing this problem. In this paper, we employ ideas similar to those in \cite{Aman00, AmWa05, shindin2026,Walk05}, specifically utilising the properties of the transport semigroup. However, in contrast to \textit{op. cit}, we focus on the spatial dependence of the diffusion coefficient and, instead of using $L_p$ spaces for the spatial part of the problem,  we build a parallel theory of the transport--fragmentation problem in $\mc X_r^0:= L_1(\mbb R_+, C_0(\Omega), dm_r)$, where $C_0(\Omega)= C(\ov{\Omega})$ if $\Omega$ is bounded, and consists of continuous functions uniformly converging to 0 at infinity if $\Omega = \mbb R^N.$  In $\mc X^0_r$, this approach requires some restrictions on the advection operator, namely, that $\mc T_0$ is only space dependent, but allows for proving the classical solvability of several classes of the transport--fragmentation--coagulation equation there.  
 
The linear section of this paper substantially generalises the results of \cite{Clov} (advection) and \cite{BanTTSP2007} (diffusion), and addresses some gaps therein. 

The paper is organised as follows. In Section \ref{MP}, we provide the theoretical backbone of the paper by showing that the structure \eqref{T0} of the transport operator allows for `gluing' semigroups generated for each $m$ in some Banach space $X_x$ along the parameter $m$ for a fairly general dependence of $\mc T_0$ on $m$, to form a semigroup in $L_1(\mbb R_+, X_x, dm_r).$  In the remainder of the paper we are concerned with $\mc X_r^1 = L_1(\mbb R_+,L_1(\Omega),dm_r)$ and $\mc X_r^0 = L_1(\mbb R_+,C_0(\Omega),dm_r).$   In Section \ref{sect3}, we first consider the transport--absorption equation, obtained by adding the loss operator $\mc A$ to $\mc T_0,$ and show the it generates a $C_0$-semigroup that has a moment regularising property. Then, in Subsection \ref{L1}, mimicking \cite{BanLam20}, we derive the existence and the moment regularising property of the transport--fragmentation in $\mc X_r^1$ for some $r$ determined by the properties of the fragmentation kernel $b$. Since some crucial $L_1$ estimates are not available in $\mc X_r^0,$ we consider an $x$ independent operator dominating the fragmentation operator in Subsection \ref{est}, where, under the assumption that the fragmentation kernel of the dominating operator is uniformly integrable with respect to the parent's mass, we show that the transport--fragmentation operator generates a positive $C_0$-semigroup in both $\mc X_r^i, i\in \{0,1\},$ under additional assumption that $\mc T_0$ is independent of $m$ if $i=0$. This result allows us to show the moment regularisation property of the transport--fragmentation semigroup also in $\mc X_r^0$. 
Finally, in Subsections \ref{sec41} and \ref{sec42} we show that the abstract results can be applied to certain classes of advection and diffusion operators, and in Subsection \ref{sec43}, we show the local classical solvability of the full transport--fragmentation--coagulation problem with unbounded coagulation kernel (controlled by the loss/absorption coefficient) if the transport operator is independent of $m$ or if the gain operator $\mc B$ is absent.


\section{Mathematical preliminaries}\label{MP}
\subsection{Spaces}
We analyse \eqref{GeneralEqn} using the theory of $C_0$-semigroups, that is, we rewrite the problem as an ordinary differential equation in appropriately chosen Banach spaces. First, let  $\Omega\subset \mbb R^N$ be either a bounded open set with a $C^2$ boundary or $\Omega = \mbb R^N$. Then, for  $i\in \{0,1\},$ we define 
\begin{subequations}
\begin{equation}
X^1_x= L_1(\Omega,dx), \label{Xxa}
\end{equation}
and
\begin{equation}
X^0_x= C_0(\Omega),\label{Xxb}
\end{equation}
\label{Xx}
\end{subequations}
where $dx$ is the Lebesque measure on $\Omega,$ and $C_0(\Omega) = C(\ov{\Omega})$ is the space of functions continuous up to the boundary $\p\Omega$ if $\Omega$ is bounded, or consists of continuous functions satisfying $u(x)\to 0$ as $\|x\|\to \infty$ if $\Omega=\mbb R^N$. Next, we introduce 
 $X_{m,r }= L_1(\mbb R_+,dm_r),$ where 
 $$
dm_r := w_r(m)dm := (1+m^r)dm,
$$
for some $r\geq 0$, where $dm$ is the Lebesque measure on $\mbb R_+$, and  the Bochner $L_1$-space
\begin{subequations}\label{mcxr}
\begin{equation}
\mc X^i_r = L_1(\mbb R_+, X^i_x, d m_r ),\label{mcxra}
\end{equation}
endowed with the norm
\begin{equation}
||u||_{\mathcal{X}^i_r}:=\int_{0}^{\infty}\|u(\cdot,m)\|_{X^i_x}dm_r.\label{mcxrb}
\end{equation}
\end{subequations}
We will omit the superscript $i$ when our considerations apply to both settings. 

In the case \eqref{Xxa}, by, e.g., \cite[Example 2.23]{BaAr} or \cite[Remark in Section 3.4]{Hille},  we have 
\begin{equation}\label{bb}
L_1(\Omega\times \mbb R_+, dx d m_r ) = L_1(\Omega, X_{m,r},dx)=L_1(\mbb R_+, X^1_x, d m_r).
\end{equation}
Alternatively, \eqref{mcxra} can be viewed as the projective tensor product, 
$\mc X^i_r = L_1(\mbb R_+)\hat \otimes_{\pi} X^i_x,
$
\cite[Section 2.3 \& Example 2.19]{Ryan}, from which, in particular, \eqref{bb} follows. Accordingly, we identify scalar functions $(t,x,m)\mapsto u(t,x,m)$ with $\mc X^i_r-$valued functions $t\to \mb u(t) := u(t,\cdot,\cdot),$ see \cite[Example 2.23]{BaAr} or \cite[Section 3.5]{Hille}. We also write $u(\cdot,m)= \mb u(m)\in X_x^i$.

We note that the spaces defined above are Banach lattices under the canonical order generated by the relation $\leq$ almost everywhere.

\begin{remark}
    The choice \eqref{Xxa} in $\mc X^1_r$ is motivated by the fact that the norms of nonnegative solutions in $L_1(\Omega\times \mbb R_+, dx dm)$ and $L_1(\Omega\times \mbb R_+, m dx dm)$ give, respectively, the total number of particles and the total mass/size of the ensemble. By considering the problem in $\mc X^i_r, r>1,$ we restrict the number of large clusters, and thus improve the properties of the fragmentation operator, see \cite{ BanDCDS20, BanLam20}. The choice \eqref{Xxb} has a mathematical motivation, and it appears in the study of nonlinear problems. 
\end{remark}
\subsection{Semigroup associated with $\mc T_0$}
Working in $\mc X_r$ defined by  \eqref{Xx} is useful in cases such as \eqref{GeneralEqn}, where $\mc T_0$ acts only on the $x$ variable with $m\in \mbb R_+$ being a parameter. Since we work in $L_1$ spaces with respect to $m$, for simplicity, if we write that a property holds for  $m\in \mbb R_+,$ we mean it holds for almost any  $m\in \mbb R_+$,  unless we want to emphasise that there are points where the statement may fail.   

If $\mc T_{0,m}$ (restricted to an appropriate domain) generates a semigroup in $X_x$ for every $m$, then \cite[Proposition 3.28]{BaAr} provides a basic criterion allowing for 'gluing' semigroups generated by the $x$ action of the operator along the parameter $m$ to obtain a semigroup in $\mc X_r.$ Here we will generalise that result in a general setting, without direct reference to the model. 

Let $X$ be a Banach space, $(\Theta,\mu)$ be a (complete) measure space and consider $\mathcal{X}=L_1(\Theta, X)$. Suppose that we are given a family of operators $\{(T_m, D(T_m))\}_{m\in\Theta}$ in $X$ and we assume that for almost every $m\in\Theta$ the operator $(T_m, D(T_m))$ generates a $C_0$-semigroup $(G_m(t))_{t\geq 0}$ in $X.$ Then, in particular,  for all $m,$ there exist constants $\omega_m $ and $M_m\geq 1$ such that
	\begin{equation}\label{G_alphaBounded}
		||G_m(t)||_{L(X)}\leq M_m e^{\omega_m t},
	\end{equation}
where $L(X)$ denotes the space of bounded linear operators on $X$.
		We define a family of operators $\mb u \mapsto \mb G(t)\mb u, t\geq 0,$ on $\mc X$ by  
	\begin{equation}\label{GtildeGen}
		[\mb {G}( t)\mb u](m)=G_m(t)\mb u(m), \quad  \mb u\in \mathcal{X},\, \quad  m\in\Theta.
	\end{equation}
	Similarly, we define 
    \begin{equation}\label{AGen}
    [\mb T\mb u](m) = T_m \mb u(m), \quad m\in\Theta,
    \end{equation}
     on 
    \begin{equation}\label{DAGen} 
    D(\mb T) := \{\mb u\in \mc X:\;\mb u(m) \in D(T_m),\,  m\in\Theta, \mb T \mb u \in \mc X\},
    \end{equation}
 and 
 \begin{equation}
[R(\la,\mb T)\mb f](m) = R(\la, T_m)\mb f(m), \quad m \in \Theta.
\label{resAa}
\end{equation}
\begin{proposition}\label{glueprop1}
 Assume that there are $M$ and $\omega$ such that for all $m$ for which $M_m$ and $\omega_m$ are defined, we have $M_m\leq M$ and $\omega_m\leq \omega$. With the above assumptions and notation,  if $\Theta\ni m \to R(\la,T_m)\mb f(m)$ is
measurable for any $\la>\omega $ and $\mb f\in \mc{X}$, then the operator
${\mb T}$ generates a semigroup $\sem{{\mb G}}$ satisfying \eqref {GtildeGen} and its resolvent is given by \eqref{resAa}. 
\end{proposition}
\noindent
\begin{proof} The proof is a modification of that of \cite[Proposition 3.28]{BaAr} that was formulated for $m$-dissipative operators $T_m$.  Here, the generation property of $T_m$ and the bounds are assumed explicitly and the fact that $\mb T$ is a generator follows from the observation that the measurability of $\mb v(m):= R(\la, T_m)\mb f(m),$ together with the uniform bound on $M_m,$ ensure that $\mb v \in \mc X$ and hence, by induction,  $m \mapsto R^n(\la, T_m)\mb f(m)$ is measurable for any $n$. Then we can integrate $\| R^n(\la, T_m)\mb f(m)\|_X,$ and the Hille-Yosida theorem and uniform estimates on $\omega_m$ gives
\begin{equation}
\|R(\la,\mb T)\mb f\|_{\mc X} = \int_0^\infty\| R^n(\la, T_m)\mb f(m)\|_X d\mu \leq \frac{M}{(\la -\omega)^n}\|\mb f\|_{\mc X},\quad \la>\omega, \mb f\in \mc X,
\label{mresrst}
\end{equation}
that is, $\mb T$ is a Hille-Yosida operator (and hence closed). To show that it is densely defined, we recall, \cite[Example 2.19]{Ryan},  that $L_1(\Theta, X)$ is isometrically isomorphic to the projective tensor product $L_1(\Theta)\widehat\otimes_\pi X$. If we consider the tensor product $L_1(\Theta)\otimes X$, 
$$
S:=\left\{\mb f = \sum\limits_{i=1}^n\phi_if_i:\; \text{for\;some\;$n$\;and\;}\phi_i\in L_1(\Theta), f_i\in X, i=1,\ldots,n\right\},$$  
 we see that $$
 S_0:=\{\mb u: \; \mb u(m)= \phi(m)u^m, \phi \in L_1(\Theta), u^m \in D(T_m)\}$$ is linearly dense in $S$. Indeed, for some $\phi f$, $\phi\in L_1(\Theta), f \in X$, we consider the Yosida approximation of $f$, $u^m_n = \la_n R\left(\la_n, T_m\right)f\in D(T_m), \la_n >\omega,$ which, by \cite[Lemma 1.3.2]{Pa}, satisfies 
 $$
 \lim\limits_{n\to \infty} u^m_n = f
 $$
 in $X$.   
 Then we have$$
 \|\phi f - \phi u^{m}_n\|_{\mc X} = \int_0^\infty \phi(m)\|f-\la_n R\left(\la_n, T_m\right)f\|_Xd\mu. 
 $$
 Since $$\|\la_n R\left(\la_n, T_m\right)f\|_X\leq \frac{\la_n}{\la_n-\omega}\|f\|_X\leq 2\|f\|_X$$
 for  $\la_n>2\omega,$ and $\phi \in L_1(\Theta)$, the Lebesgue dominated convergence theorem gives
 $$
 \lim\limits_{n\to \infty}\|\phi f - \phi u^{m}_n\|_{\mc X} = 0.
 $$
 Next, by linearity, 
 $$
 \sum\limits_{i=1}^n \phi_i(m)u^{m}_{i,n}\in D(\mb T),
 $$
 for $\phi_i \in L_1(\Theta)$, so $D(\mb T)$ is dense in $S$. Since, as in the aforementioned \cite[Example 2.19]{Ryan}, $S$ is dense in $\mc X,$  $D(\mb T)$ is also dense in $\mc X$. \end{proof}

As with the advection semigroup, the semigroup's formula is simpler than its resolvent's. Thus, we state and prove a semigroup version of the above proposition. 
     \begin{proposition}\label{ParameterSemigroup}
	If for any $\mb u \in \mc X$ and any fixed $t\geq 0,$ the  function $m\mapsto G_m(t)u(m)$ is 
    Bochner-measurable and there exists a function  $\Theta\ni t\mapsto M(t)\in \mbb R_+$ such that for any $t\geq 0$ $\sup_{m \in \Theta} M_m e^{\omega_m t}\leq M(t)$,  then $(\mb {G}(t))_{t\geq 0},$ given by \eqref{GtildeGen}, is a $C_0$-semigroup generated by $\mb T$.
\end{proposition}
\begin{proof}
	For any $t\geq0$ and $\mb u \in\mathcal{X},$ the function  
    $m\mapsto G_m(t)\mb u(m)$ is Bochner-measurable, hence, by, e.g., \cite[Theorem 2.22]{BaAr}, $m\mapsto \|G_m(t)u(m)\|_X$ is Lebesgue measurable. Hence, 
    by \eqref{G_alphaBounded}, and, again, \cite[Theorem 2.22]{BaAr}, $m\mapsto G_m(t)\mb u(m)$ is Bochner-integrable giving, 
    $$
    \|\mb G(t)\mb u\|_{\mc X}\leq M(t)\|\mb u\|_{\mc X}, \quad t\geq 0.
    $$
    So,  
	 $\mb {G}(t)\in\mathcal{L}(\mathcal{X})$ for $t\geq 0$, that is, $(\mb {G}(t))_{t\geq 0}$ is a family of bounded linear operators on $\mathcal{X}$.

The algebraic semigroup properties are easily verified. To prove strong continuity, we get
			$$
		||G_m(t)\mb u(m)-\mb u(m)||_X\leq \|G_m(t)\mb u(m)\|_X + \|\mb u(m)\|_X \leq (M(t)+1) \|\mb u(m)\|_X =: g(m),
			$$
			where $g \in L^1(\Theta)$ since $\pmb u \in \mathcal{X}$. Thus, the strong continuity of $(\mb {G}(t))_{t \geq 0}$ follows since $(G_m(t))_{t\geq 0}$ are $C_0$-semigroups for $m\in \Theta$, and by the Lebesgue Dominated Convergence theorem.
					
To prove the last statement, we use \cite[Remark 2.23, Eqn (2.36)]{BaAr} (or an extension of \cite[Eqn (3.4.6)]{Hille} to an infinite interval) to write the resolvent equation
\begin{equation}
[R(\la,\mb T)\mb u](m) = \left[\int_0^\infty e^{\la t}\mb G(t)\mb u dt \right](m) = \int_0^\infty e^{\la t} G_m(t)\mb u(m) dt = R(\la, T_m)\mb u(m), 
\label{resm}
\end{equation}
valid for $m \in \Theta$ and sufficiently large $\la\in \mbb R,$ to claim the equality of generators. 
\end{proof} 
Using \eqref{mresrst} and \eqref{resm}, we immediately have 
\begin{corollary}
In cases of Propositions \ref{glueprop1} and \ref{ParameterSemigroup}, if $(G_m(t))_{t\geq 0}$ is analytic, then so is $(\mb {G}(t))_{t \geq 0}$.\label{ana1}
\end{corollary}
	Next, we provide a practical criterion for the measurability of the resolvent. 
\begin{proposition}\label{propmes}
Let $\Theta$ be an open set and consider a family of operators $((T_m,D))_{m\in \Theta}$ on a common domain $D$. Let, for any $m \in \Theta,$ there exist a sequence $((T_{m,n},D))_{n \in \mbb N}$ of  
operators such  that for each $u\in D$ and $n\in \mbb N$, $m\mapsto T_{m,n}u$  is a continuous function and $\lim_{n\to \infty}T_{m,n}u = T_m u$. Assume that the resolvents of $T_m, T_{m,n}, m \in \Theta,n\in \mbb N,$ are defined for $\la>\la_0$ for some $\la_0,$ and are uniformly bounded: there exists $C$ such that for almost all $m \in \Theta, n \in \mbb N,$ and $\la >\la_0$ 
\begin{equation}
\|R(\la, T_m)\|_{L(X)}\leq C, \quad \|R(\la, T_{m,n})\|_{L(X)}\leq C.
\label{resest1}
\end{equation}
Then, for any $\mb f\in \mc X$,  
$$
\Theta\ni m\mapsto R(\la, T_m)\mb f(m)\in X
$$
is Bochner-measurable. 
\end{proposition}
    \begin{proof} First, let us consider a family $((T_m,D))_{m\in \Theta}$ of operators such that $m\mapsto T_m u$ is continuous for every $u\in D$, and let $m\mapsto \mb f(m)$ be a continuous $X$-valued function.  For a given $m\in \Theta$ and $h$ such that $m+h\in \Theta$ consider 
  \begin{align*}
  \la  \mb u(m) - T_m \mb u(m) &= \mb f(m),\\
  \la \mb u(m+h) - T_{m+h} \mb u(m+h) &= \mb f(m+h).
  \end{align*}
  Subtracting, rearranging (made possible by the common domain $D$), and applying the resolvent, we get 
     \begin{equation}
     \mb u(m+h)-\mb u(m) = R(\la, T_{m+h})(T_{m+h}-T_m)\mb u(m) + R(\la, T_{m+h})(\mb f(m+h)-\mb f(m)).
     \label{uRf1}
     \end{equation}
     Taking the norm and using \eqref{resest1} gives 
     \begin{align*}
     \|\mb u(m+h)-\mb u(m)\|_{X} &\!\leq \!\| R(\la, T_{m+h})(T_{m+h}-T_m)\mb u(m)\|_{X} \!+ \|R(\la, T_{m+h})( \mb f(m+h)-\mb f(m))\|_{X}\\
     &\!\leq\! C(\|(T_{m+h}-T_m)\mb u(m)\|_{X} + \|\mb f(m+h)-\mb f(m)\|_{X}).
     \end{align*}
     The continuity of $m\mapsto \mb f(m)$ and $T_m$ implies that $m \mapsto R(\la, T_m)\mb f(m)$ is continuous (and hence Bochner-measurable). Next, we remove the assumption of the continuity of $\mb f$. Let us consider $\mb f\in \mc X$. As in the proof of Proposition \ref{glueprop1} (or using directly \cite[p.13]{ABHN}), for any $\mb f \in \mc X$ there is a sequence $(\mb f_n)_{n\in N}$ of  continuous (even $C_0^\infty$) functions with respect to $m$ such that
     $$
     \lim\limits_{n\to\infty}\int\limits_\Theta \|\mb f_n(m)-\mb f(m)\|_Xd\mu = 0.
     $$ 
     Since the integrands are real functions, there is a subsequence of $(\mb f_n)_{n\in \mbb N}$ converging to $\mb f$ almost everywhere in $X$. Passing to this subsequence without changing indices,  for $m \in \Theta$ we have 
     \begin{equation}
     \lim\limits_{n\to \infty}\mb f_n(m) = \mb f(m)\quad \text{and}\quad \lim\limits_{n\to \infty}R(\la, T_m)\mb f_n(m) = R(\la, T_m)\mb f(m)
     \label{uRf}
     \end{equation}
     in $X$.
     On the other hand, 
     $$
     \mb u_n(m) = R(\la,T_m)\mb f_n(m), \quad n\in \mbb N,
     $$
     are continuous functions and hence we can evaluate the norm $\|\mb u_n-\mb u_m\|_{\mc X},$ which gives  
     $$\|\mb u_n-\mb u_m\|_{\mc X} \leq C\|\mb f_n-\mb f_m\|_{\mc X}.$$
     Thus, there is $\mb u\in \mc X$ such that 
     $$
     \lim_{n\to \infty}\int_0^\infty\|\mb u(m)-\mb u_n(m)\|_X d \mu =\lim_{n\to \infty}\int_0^\infty\|\mb u(m)-R(\la, T_m)\mb f_n(m)\|_X d \mu =0.$$  
     Hence, as above, there is a subsequence $(n_k)_{k\in \mbb N}$ such that for $m\in \Theta$
     $$
     \mb u(m) = \lim\limits_{k\to \infty} R(\la, T_m)\mb f_{n_k}(m)
     $$
     in $X$. Using \eqref{uRf}, we see that
     $$
     \mb u(m) = R(\la, T_m)\mb f(m),
     $$
     thus, $m\mapsto R(\la, T_m)\mb f(m)\in \mc X$ for any $\mb f \in \mc X,$ (and so, in particular, it is Bochner-measurable). Finally, we drop the assumption of the continuity of $T_m$ and consider 
     \begin{align*}
  \la \mb u(m) - T_m \mb u(m) &= \mb f(m),\\
  \la \mb u_n(m) - T_{m,n} \mb u_n(m) &= \mb f(m), 
  \end{align*}
  where $T_{m,n}$ are as in the assumptions of the proposition. Thus, as in \eqref{uRf1}, 
   $$
  \mb u_n(m) - \mb u(m) =R(\la, T_{m,n}) (T_{m,n}- T_{m})\mb u(m).
  $$  
  Taking norms, we get 
  $$
  \|\mb u_n(m)-\mb u(m)\|_{X}\leq C\|(T_{m,n}-T_{m})\mb u(m)\|_{X},
  $$
  hence, by assumption, for almost any $m$, $\mb u(m) = \lim_{n\to \infty}\mb u_n(m)$ in $X$. Since $\mb u_n, n\in \mbb N,$ are Bochner-measurable by the first part of the proof, $\mb u$ is also Bochner-measurable.\end{proof}
We observe that the assumption of the continuity of $m\mapsto T_{m,n}u$ was necessary only in the first part of the proof and can be skipped as long as we know that 
$m\mapsto R(\la,T_{m,n})\mb f(m)$ are measurable. Then we can repeat the estimates from the last part of the proof to get 
\begin{corollary}\label{cormes}
Assume that the operators  $T_{m,n}, {n \in \mbb N}, m \in \Theta,$ satisfy all assumptions of Proposition \ref{propmes} except for the continuity in $m$ but are such that   
  $m \mapsto \mb u_n(m): =\! R(\la,T_{m,n})\mb f(m)$ are Bochner-measurable for every $\mb f\in \mc X$. 
Then, 
$$
\Theta\ni m\mapsto R(\la, T_m)\mb f(m)\in X
$$
is Bochner-measurable for any $\mb f\in \mc X$. \end{corollary}

\section{Solvability of the transport--fragmentation model}\label{sect3}
	To prove the well-posedness of the transport--fragmentation problem, we proceed in two steps. First, we prove the solvability of the transport--advection problem using the results of the previous section, and then we apply the Desch--Voigt perturbation theorem to show the well-posedness of the transport--fragmentation semigroup. 
    
	\subsection{The transport--absorption  problem}\label{sect31}
	We begin our analysis by considering 
		\begin{equation}\label{TransProblem}
		\begin{split}
	\partial_t u(t,x,m)&= \mc Tu(t,{x},m):=\mc T_0u(t,{x},m) +\mc Au(t,{x},m),\\
			u(0,x,m)&=\mr{u}(x,m), \hspace{0.5cm} \text{a.e.}~ (x,m)\in\Omega\times\mathbb{R}_+,
		\end{split}
	\end{equation}
	    see \eqref{ops2}. As noted above, here we treat $m$ as a parameter and thus work in ${X}^i_x$. We assume 
\begin{subequations}\label{Assumptionofa}
\begin{equation}
0\leq a \in L_{\infty,loc}(\Omega\times\mathbb{R}_+),
\label{aa0}
\end{equation}
in the case \eqref{Xxa}, 
\begin{equation}
0\leq a \in L_{\infty,loc}(\mbb R_+, C_b(\Omega)),
\label{aa0c}
\end{equation}
where $C_b(\Omega)$ is the set of bounded continuous functions on $\Omega$, in the case \eqref{Xxb},
and, since we need some uniformity of $a$ with respect to  $x,$ we assume that there are $\alpha_1, \alpha_2\in L_{\infty,loc}(\mathbb{R}_+)$ and $0<M<\infty$ such  that 
\begin{equation}
\alpha_1(m)\leq a({x},m)\leq \alpha_2(m)\leq M\alpha_1(m), \quad m \in \mbb R_+. \label{aa1}
\end{equation}
\end{subequations}
We shall often need to control the convergence of $e^{-\alpha_1(m) t}$ to 0 as $t\to \infty$. Thus, often we additionally  assume  
that there exist positive constants $a_0,~\gamma,$ and $m_0\geq 1$ such that 
				\begin{equation}\label{aCondition}
					\alpha_1(m)\geq a_0m^{\gamma}, \hspace{0.5cm}  m\geq m_0.
				\end{equation}
    We observe that $\mc T$ is also a family of operators $(\mc T_m)_{m\in \mbb R_+}$ on $X_x$, satisfying \eqref{T0}. 
    
     For $m \in \mbb R_+,$  we consider the restriction $T^i_{0,m}$ of $\mc T_{0,m}$, see \eqref{T0}, to a domain $D(T^i_{0,m})\subset X^i_x,$ and let $[A^i_m u](x,m) :=-a(x,m)u(x,m)$; by \eqref{aa1}, $A^i_m$ is defined on $X^i_x$ for any $m$. 
     
     To formulate the next results, let us recall that a positive $C_0$-semigroup on a Banach lattice is called substochastic if it is contractive, and stochastic if it is conservative on the positive cone. Using these definitions, we assume that $(T^i_{0,m}, D(T^i_{0,m}))$ generates a substochastic semigroup, say $(G_{T^i_{0,m}}(t))_{t\geq 0},$ in $X^i_x.$ Then, under assumption \eqref{aa0}, respectively, \eqref{aa0c}, the semigroup $(G_{T^i_{m}}(t))_{t\geq 0}$, generated by  $(T^i_m, D(T^i_{0,m})) = (T^i_{0,m} +A^i_m, D(T^i_{0,m})),$ is also substochastic in respective $X^i_x$. 
     Further, let 
     \begin{equation}\label{dom}
     \begin{split}
     D(\mb T^i_0) &= \{\mb u\in \mc X^i_r:\; \mb u(m)\in D(T^i_{0,m}), \mc T_0u \in \mc X^i_r\}, \\
     D(\mb A^i) &= \{\mb u\in \mc X^i_r:\; au \in \mc X^i_r\},\\
     D(\mb T^i)&=D(\mb T^i_0)\cap D(\mb A^i),
     \end{split}
     \end{equation}
     with $\mb T^i_0=\mc T_0|_{D(\mb T^i_0)}, \mb A^i=\mc A|_{D(\mb A^i)}$ and $\mb T^i=\mc T|_{D(\mb T^i)}$. 
     
     We note that the above operators are realisations of the corresponding expressions in $\mc X_r$, so they depend on $ r$. However, since we are working with a fixed $r$ for the time being, we will ignore this dependence in the notation unless it is relevant.

     The main assumption is 
     \begin{description}
     \item {(A1)} $(\mb T^i_0,D(\mb T^i_0))$ generates a substochastic semigroup in $\mc X^i_r$, say $(G_{\mb T^i_0}(t))_{t\geq 0},$ satisfying 
     $$
     [G_{\mb T^i_0}(t)\mb u](x,m) = [G_{T^i_{0,m}}(t)\mb u(m)](x), \quad m\in \mbb R_+.
     $$
     \end{description}
     \begin{proposition}\label{prop31}
     If assumption (A1) is satisfied and $\mb A^i$ is defined as above with $a$ satisfying the respective version of \eqref{Assumptionofa}, then  $(\mb T^i,D(\mb T^i))$ generates a substochastic semigroup in $\mc X^i_r$, say $(G_{\mb T^i}(t))_{t\geq 0},$ satisfying 
     $$
     [G_{\mb T^i}(t)\mb u](x,m) = [G_{T^i_{m}}(t)\mb u(m)](x),\quad m\in \mbb R_+.
     $$
          \end{proposition}
          \begin{proof}
          Let $a$ satisfy \eqref{Assumptionofa} with \eqref{aa0} or \eqref{aa0c} determined by the choice of $X^i_x.$  Hereafter, we skip $i$ in the notation unless it is relevant. Define $a_n (x,m)= \chi_{[0,n]}(m)a(x,m)$, where $\chi$ is the characteristic function of $[0,n]$. Then $0\leq a_n(x,m)\leq \sup_{0\leq m\leq n}\alpha_2(m)<\infty.$ 
          Hence, $a_n$ is a bounded measurable function on $\Omega\times \mbb R_+$, and,  if  \eqref{aa0c} is satisfied,  it is continuous in $x$.    Then $\mb T_0+\mb A_n$, where $\mb A_n$ is the operator of multiplication by $-a_n,$ generates a positive semigroup of contractions on  $\mc X_r$. In particular, $m\mapsto R(\la, \mb A_n)\mb f(m)$ is a measurable $X_x$-valued function for any $\mb f\in \mc X_r$. and,  proceeding as in the proof of Proposition \ref{propmes}, for $\mb f\in \mc X_r$, and almost any $(x,m)\in \Omega\times \mbb R_+$, 
          \begin{equation}
          \begin{split}
  \la  u(x, m) - [T_{0,m}  u(\cdot, m)](x) + a(x,m) u(x,m)&= f(x,m),\\
  \la  u_n(x,m) - [T_{0,m} u_n(\cdot,m)](x)+ a_n(x,m)u_n(x,m) &= f(x,m).
  \end{split}
  \label{AR0}
  \end{equation}
    Since $D(T_{0,m}+A_{m,n})$ does not depend on $n$, subtracting, re-arranging and applying the resolvent, we get
         $$
  u(x,m) - u_n(x,m) =[R(\la, T_{0,m}+A_{m,n})((a(\cdot,m)-a_n(\cdot,m))u(\cdot,m))](x),
  $$  
  and, using the fact that $T_{0,m}+A_{m,n},m\in \mbb R_+,$ generate contraction semigroups, 
  \begin{equation}
  \| u(\cdot,m)- u_n(\cdot,m)\|_{X_x}\leq \|(a(\cdot,m)-a_n(\cdot,m))u(\cdot,m)\|_{X_x}.
  \label{uun}
  \end{equation}
  Now, for a given $m,$  
  \begin{equation}
  a(x, m) - a_n(x, m) = \chi_{[n,\infty)}(m)a(x,m) =0, \quad n>m,\label{aan}
  \end{equation}
  for any $x\in \Omega.$ In case \eqref{aa0c}, this means that for a. a. $m$
  $$
  \lim\limits_{n\to\infty} \|(a(\cdot,m)-a_n(\cdot,m))u(\cdot,m)\|_{X_x} = 0,$$
  and hence \eqref{uun} implies the thesis. 
  
  In the case \eqref{aa0}, \eqref{aan} implies that for $m\in \mbb R_+$, $\lim_{n\to \infty}a_n(x,m)u(x,m) = a(x,m)u(x,m)$ almost everywhere on $\Omega$ and $ |(a(\cdot,m)-a_n(\cdot,m))u(\cdot,m)|\leq 2\alpha_2(m)u(\cdot,m) \in L_1(\Omega,dx)$. Hence, by the Dominated Convergence Theorem, $\lim_{n\to \infty}u_n(\cdot,m) = u(\cdot,m)$ in $X_x$ for almost any $m.$ Thus, in both cases,  since $\mb u_n$ are Bochner measurable, so is $\mb u$.  
  
  Hence,  Proposition \ref{glueprop1} implies the existence of a semigroup generated by  $\mbb T = (\mc T_{0,m}+\mc A)|_{D(\mbb T)}$ on  
  $D(\mbb T) = \{\mb u\in \mc X_r:\; \mb u(m) \in D(T_{0,m}), (x,m)\mapsto \mc T_{0,m} u(x,m)-a(x,m)u(x,m)\in \mc X_r\}$, see \eqref{DAGen}. To show that $\mbb T=\mb T$, we note that since 
  $$
  T_{0,m} - \alpha_1(m) I  = T_m +(a(x,m)-\alpha_1(m)) I,
  $$
 $T_m$ is a generator and $a(x,m)-\alpha_1(m)\geq 0$, 
  we obtain, as in the Bounded Perturbation Theorem,    
  \begin{equation}
     R(\la, T_{m})\leq R(\la, T_{0,m}- \alpha_1(m) I), \quad m\in \mbb R_+.\label{Tmalpha}
  \end{equation}
 Thus, using the fact that $T_{0,m}, m\in \mbb R_+,$ are dissipative on $X_x$,  for any $\mb f\in \mc X_r$, 
 \begin{equation}
 \begin{split}
  \| A_m   R(\la, T_{m})\mb f(m)\|_{X_x}&\leq M\|\alpha_1(m)R(\la, T_{0,m}- \alpha_1(m) I)\mb f(m)\|_{X_x}\\
  &\leq \frac{M\alpha_1(m)}{\la+\alpha_1(m)}\|\mb f(m)\|_{X_x}\leq M\|\mb f(m)\|_{X_x}.
  \end{split}\label{alph1est}
  \end{equation}
  From the first part of the proof, we know that $m\to \mb u(m):= [R(\la, \mbb T)\mb f](m)$ is Bochner-measurable and thus, upon integration, we see that $\mb u \in \mc X_r$. Now, if $i=1$, then, from \eqref{bb}, $(x,m)\mapsto u(x,m)$ is measurable and hence $ (x,m)\mapsto a(x,m)u(x,m)$ is also measurable, and, by \eqref{aa0}, $au \in L_1(\Omega\times K, dxdm_r)$ for any bounded $K\subset \mbb R_+$. Thus, $m\mapsto a(\cdot,m)u(\cdot,m)$ is $X_x^1$-measurable. 
  For $i=0$, we see that if $a_k\in C_b(\Omega)$ and $u_j\in C_0(\Omega)$,  then $a_ku_j\in C_0(\Omega)$. Hence, by \eqref{aa0c} and the definition of Bochner measurability, $m\mapsto a(\cdot,m)u_j$ is $X_x^0$-measurable and then, since, for $m\in \mbb R_+$ 
\begin{equation}
a(\cdot,m)u(\cdot,m) = a(\cdot,m)\lim\limits_{n\to\infty}\sum\limits_{j=1}^n\chi_{I_{nj}}(m) u_{nj} = \lim\limits_{n\to\infty}\sum\limits_{j=1}^n\chi_{I_{nj}}(m) a(\cdot,m)u_{nj},
\label{mesau}
\end{equation}
  where $\chi_{I}$ is the characteristic function of interval $I,$ and $u_{nj}\in X_x^0,$ $au$ is also $X_x^0$-measurable as an a.e. limit of measurable functions.
  Therefore, in both cases, we can integrate \eqref{alph1est} with respect to $dm_r$, getting 
  \begin{equation}
  \|\mb A R(\la, \mbb T)\mb f\|_{\mc X_r} \leq M\|\mb f\|_{\mc X_r},
  \label{AR1}
  \end{equation}
 and hence $\mb u \in D(\mb A)$. Using the first equation of \eqref{AR0}, we see that $\mb u \in D(\mb T_0)$. Thus, $D(\mbb T)\in D(\mb A)\cap D(\mb T_0)$. Since the reverse inclusion is obvious, we obtain $\mbb T=\mb T$.
          \end{proof}
Since $a(\cdot,m)$ is bounded for $m\in \mbb R_+$, if $(G_{T_{0,m}}(t))_{t\geq 0}$   is analytic, then so is $(G_{T_{m}}(t))_{t\geq 0}$, hence, by Corollary \ref{ana1}, we have
\begin{corollary}
If the semigroup $(G_{T_{0,m}}(t))_{t\geq 0}$ is analytic for any $m\in \mbb R_+$, then so is $(G_{\mb T}(t))_{t\geq 0}$. \label{ana2}
\end{corollary}
     
     In the case \eqref{Xxa}, we have the following improvement of \eqref{AR1}. 
     \begin{lemma}\label{lemAR} Let $X_x = X_x^1 = L_1(\Omega,dx)$ and \eqref{aa0} be satisfied. Then,
     \begin{equation}\label{AuEst2}
||\mb A R(\lambda,\mb T)\mb f||_{\mc X^1_r}\leq||\mb f||_{\mc X^1_r}, \quad \mb f\in \mc X^1_r.
				\end{equation}
     \end{lemma}
     \noindent
       \begin{proof} Since $(G_{T^1_{0,m}}(t))_{t\geq 0}, m\in \mbb R_+,$ are positive semigroups of contractions, we have 
        \begin{equation}\int\limits_0^\infty\int\limits_{\Omega}\mb T_{0}\mb u dx dm_r = \int\limits_0^\infty\left(\int\limits_{\Omega}T^1_{0,m}\mb u(m) dx\right) dm_r \leq 0\label{T0dis}
        \end{equation} 
        for any $u \in D(\mb T^1)_+.$         
        Since $u(\cdot,\cdot)=\mb u,$ where $D(\mb T^1)_+ \ni \mb u = R(\la, T^1)\mb f$, $\mb f\in \mc X^1_{r,+},$ satisfies term-wise the resolvent equation
        \begin{equation}
  \la u(x,m) - [\mc T_0 u](x,m) +a(x,m) u(x,m) = f(x,m),
    \label{req01}
  \end{equation}
        integrating and using the fact that $T_{0,m}$ is dissipative,  we obtain
\begin{align*}\int\limits_{0}^{\infty} \int\limits_{\Omega} (\lambda u - \mc T_0u + au )dx dm_r  &= -\int\limits_{0}^{\infty} \int\limits_{\Omega} \mb T^1_0\mb u dx dm_r  	 + \int\limits_{0}^{\infty} \int\limits_{\Omega} (\lambda u + au )dx dm_r = \int\limits_{\Omega}\int\limits_{0}^{\infty} f dx dm_r,\end{align*} that is, by $-\int_{0}^{\infty} \int_{\Omega} \mb T^1_0\mb u dx dm_r\geq 0$, 
\begin{equation}
||\mb A^1\mb u||_{\mc X^1_r} + \lambda ||\mb u||_{\mc X_r}	\leq  ||\mb f||_{\mc X^1_r}.
\label{AR}
\end{equation}
Since $\mb u=R(\lambda,\mb T)\mb f,$ extending \eqref{AR} to $\mc X^1_r=\mc X^1_{r,+}-\mc X^1_{r,+}$, we get  \eqref{AuEst2}. \end{proof}
Moreover, in the general case, we get 
     \begin{proposition} \label{lemGTmest} If $a$ satisfies \eqref{Assumptionofa}, then for almost all $m\in \mbb R_+,$
     \begin{equation}
     \|G_{T_m}(t) f\|_{X_x}\leq e^{-t\alpha_1(m)} \|f\|_{X_x}, \quad f\in X_x. 
     \label{GTmest0}
     \end{equation}
     Let, additionally, \eqref{aCondition} be satisfied and  $r:=p+q, p\geq 0.$ For any $q\geq 0$ there exist constants $C_1,C_2$ such that  
     \begin{equation}
     \| G_{\mb T}(t) \mb f\|_{\mc X_{r}}\leq \left( C_1+\frac{C_2}{t^{\frac{q}{\gamma}}}\right)\|\mb f\|_{\mc X_p}, \quad f\in \mc X_p. 
     \label{GTmest}
     \end{equation}
      \end{proposition}
     \begin{proof}
   Using \eqref{Tmalpha} and the exponential formula for semigroups, \cite[Corollary III 5.5]{EN}, we find that for $f\in X_{x,+}, m \in \mbb R_+, t\geq 0,$
   $$0\leq G_{T_m}(t)f\leq G_{T_{0,m}-\alpha_1(m) I}(t)f = e^{-\alpha_1(m) t}G_{T_{0,m}}(t)f,$$
   from which \eqref{GTmest0} follows by the contractivity of $\sem{G_{T_{0,m}}}$ in $X_x$.

Now, we prove \eqref{GTmest}. Let $\mb f \in \mc X_p$. Recalling that $m_0\geq 1$ in \eqref{aCondition}, and $r=p+q,$ 
\begin{equation*}
\begin{split}
\| G_{\mb T}(t)\mb f\|_{\mc X_{r}} &\leq 
\int\limits_{0}^{m_0} e^{-\alpha_1(m)t}\|f(\cdot,m)\|_{X_x}(1+m^{q}) dm_p
+\int\limits_{m_0}^{\infty} e^{-a_0m^{\gamma}t}\|f(\cdot,m)\|_{X_x}(1+m^{q}) dm_p\\ 
&\leq  (C_1 + 2\max\limits_{m\in \mbb R_+}e^{-a_0m^{\gamma}t}m^{q})\|\mb f\|_{\mc X_p},
\end{split}
\end{equation*}
which, upon evaluating 
$$
\max\limits_{m\in \mbb R_+} e^{-a_0m^{\gamma}t}m^{q} =\frac{1}{t^\frac{q}{\gamma}}\max\limits_{z\in \mbb R_+}   e^{-a_0z}{z}^{\frac{q}{\gamma}} =\frac{1}{t^\frac{q}{\gamma}} e^{-\frac{q}{\gamma}}\left(\frac{q}{\gamma a_0}\right)^\frac{q}{\gamma} =:\frac{C_2}{2}\frac{1}{t^\frac{q}{\gamma}},
$$
gives \eqref{GTmest}. 
     \end{proof}
          
\subsection{Full transport--fragmentation equation} \label{sectTF}
   In this section, we analyse the transport--fragmentation part of \eqref{GeneralEqn},  \begin{equation}\label{ACP1}
	\begin{split}
	\partial_t u(t,{x},m)&=\mc T_0u(t,{x},m)-a({x},m)u(t,{x},m)+\int_{m}^{\infty}b({x},m,s)a({x},s)u(t,{x},s)ds,\\
		u(0,{x},m)&=\mr{u}({x},m),
	\end{split}
\end{equation}
for $t>0,({x},m)\in\Omega\times\mathbb{R}_+$. For this, we need to discuss the gain operator, defined by  $$\mc Bu(x,m) = \int_m^\infty a(x,m)b(x,m,s)u(x,s)ds$$ restricted to a suitable domain. 
We begin with the properties of the fragmentation kernel $b$. We assume that  $b\geq0 $ is a measurable function satisfying $b(x,m,s)=0$ for $m>s, x\in \Omega$. For each $r\geq 0$ we define the $r$th moment of $b$ and its deviation from $s^r$ by, respectively,
\begin{subequations}
		\begin{equation}\label{defnOfnTilde}
			n_r({x},s):=\int_{0}^{s}b({x},m,s)m^rdm, 	
		\end{equation} 
\begin{equation}
 N_r({x},s):=s^r-n_r({x},s).
\label{Nr}	
\end{equation}
\end{subequations}
The total mass of the daughter particles is given by $n_1,$ and if we assume that no mass is lost or created in the process, we must have 
  \begin{equation}\int_{0}^{s}mb(x,m,s)dm=s, \quad x\in \Omega.\label{b1}
    \end{equation}
The expected number of daughter particles produced by the fragmentation of a mass $s$ particle is given by $n_0.$ We assume that there are constants $b_0\geq 1$ and $l\geq 0$ such that for a.a. $x\in \Omega$, 
		\begin{equation}\label{n(s)bounded}
			n_0({x},s)\leq b_0(1+s^l).
		\end{equation}
The fact that $b_0\geq 1$ follows from $n_0(x,s)\geq 1,$ implied by \eqref{b1}.  
        
Now, we split the considerations into two streams, dealing separately with $X_x^1$ and $X_x^0.$
\subsection{$L_1$ theory}\label{L1}

Thanks to \eqref{bb}, the problem in $\mc X^1_r$ does not significantly differs from the space homogeneous one, as we can interchange the order of taking norms in $X_x^1$ and $X_{m,r}$ and hence the proofs of the main theorems of this section are almost identical to the proofs of analogous results in \cite{BanLam20} and thus will be omitted.  

Standard calculations, see \cite[Section 5.1.7]{BLL1}, show that $\mc B$ restricted to $D(\mb T)$ defines a positive operator in $\mc X^1_r$, which we denote by $\mb B$.

We take an arbitrary $r$ such that  
\begin{equation}\label{r>max1l}
				r>\max\{1,l\}, 	
			\end{equation}
see \eqref{n(s)bounded}, and introduce the space uniform version of the basic assumption, allowing for the proof that $\mb B$ is the Miyadera--Desch perturbation of $\mb T,$ see 
  \cite[Sections 5.1.7 and 5.2.3]{BLL1}, \cite{BanDCDS20, BanLam20}, that is, we assume that for $r$ satisfying \eqref{r>max1l} there exist $c_r<1$ and $s_r>0$ such that
            \begin{equation}\label{c_r(x)Bound}
	n_r(x,s)\leq c_r s^r, \quad  s\geq s_r.
			\end{equation}
            Then, thanks to \eqref{AuEst2}, the following theorem can be proved exactly as \cite[Theorem 2.2]{BanLam20}. 
\begin{theorem}\label{FGeneratesSemigroup}
Let \eqref{Assumptionofa}, \eqref{n(s)bounded}, \eqref{r>max1l}, and  assumption (A1) 
be satisfied.	Then $(\mb K,D(\mb T)):=(\mb T+\mb B, D(\mb T))= (\mb T_0+\mb A+\mb B,D(\mb T))$ generates a positive $C_0$-semigroup, $( G_{\mb K}(t))_{t\geq 0}$, on $\mathcal{X}^1_r$.
\end{theorem}
 Using again the fact that we can interchange the order of integration when we derive inequalities for the moments of the solution, 
            $$
            M_r(t) = \int\limits_{\Omega}\int\limits_{0}^\infty u(t,x,m) dm_r dx,
            $$
            we can first integrate the RHS of \eqref{ACP1} with respect to $x$ to eliminate the contribution of the differential operator due to its dissipativity, and then change the order of integration to proceed with the evaluation of the integrals of the remaining terms, as in \cite[Lemma 5.1.34 \& Theorem 5.1.48]{BLL1}.  This leads to 
            \begin{equation}\label{d/dtM0r}
					\begin{split}
						\dfrac{d}{dt}M_{r}(t)&\leq -\int_{\Omega}\int_{0}^{\infty} (N_0(x,m)+N_r(x,m))a(x,m)u(t,x,m)dmdx, \quad \mb u \in D(\mb K).	
					\end{split}
				\end{equation}
For the next result, we will need the scale of spaces $\mc X^1_r, r\geq 0.$  By \cite[Proposition 5.1.33]{BLL1}, $\mathcal{X}^1_{r_2}$ is continuously embedded in $\mathcal{X}^1_{r_1}$ if $r_1<r_2$. We will slightly abuse the notation and use the same symbols for operators in $\mc X_r^1$ with different $r$, but will distinguish them by domains, e.g.,   
\begin{equation}
D_r(\mb K) = D_r(\mb T)= \{\mb u\in \mc X_r^1:\; \mc T_0\mb u\in \mc X^1_r\;\text{and}\;\mc A u \in \mc X^1_r\}.\label{Dr}
\end{equation}
The semigroup $\sem{G_{\mb K}}$ operating in $\mc X_{r_2}^1$ is the restriction of $\sem{G_{\mb K}}$ in $\mc X^1_{r_1}.$
				Inequality \eqref{d/dtM0r} is exactly  \cite[Eqn (2.35)]{BanLam20} (without the growth coefficient),  so that we can repeat the proof of \cite[Theorem 2.3]{BanLam20} to establish the following theorem.     		        
            \begin{theorem}\label{G_Kbounded}
				Let the assumptions of Theorem  \ref{FGeneratesSemigroup} hold and assume that $a$ satisfies \eqref{aCondition}. 
				Then, for any $n,r$ and $q$ satisfying $\max\{1,l\}<n<p<r,$ there are $C>0$ and $\theta>0$  such that 
	\begin{equation}
	|| G_{\mb K}(t)\mb {\mr u}||_{\mc X^1_r}\leq Ce^{\theta t}t^{\frac{n-r}{\gamma}}||\mb{\mr u}||_{\mc X^1_p}, \hspace{0.5cm} \text{ for all } \mb {\mr u}\in\mathcal{X}^1_p.
\label{momreg}				
\end{equation}
			\end{theorem}

\subsection{General $X_x$ theory}\label{est}

In this section, we consider the $X_x^0$ theory of fragmentation. The approach developed here can also be used in $X_x^1$, and we will apply it to improve \eqref{momreg} if $\mc T_0$ is independent of $m$. 

We assume \eqref{Assumptionofa}, and let $\beta \geq 0$ with $\text{supp}\, \beta (m,s)\subseteq \Delta:=\{(m,s)\in {\mathbb R}^2_+;\; m\leq s\}$ be a measurable function. 
Then, we consider 
\begin{equation}\label{TransProblem2}
		\begin{split}
	\partial_t u(t,x,m)&=\mc T_0 u(t,x,m) -\alpha_1(m)u(t,x,m) + \int_m^\infty \alpha_2(s)\beta(m,s)u(t,x,m)ds\\
    &=: \mc T_0 u(t,x,m) + \mc A_1 u(t,x,m) + \mc B_1 u(t,x,m),\quad t>0,\, (x,m)\in\Omega\times\mathbb{R}_+, \\
			u(0,x,m)&=\mr{u}(x,m), \quad (x,m)\in\Omega\times\mathbb{R}_+. 
		\end{split}
	\end{equation}
    We denote by ${\sf A}$ the operator of multiplication by $\alpha_1$ defined on $D({\sf A}) = \{u\in X_{m,r}:\;\alpha_1 u \in X_{m,r}\},$ and by ${\sf B}$ the restriction of $\mc B_1$ to $D({\sf A}),$ which, as we prove in Proposition \ref{propei1}, is well defined for sufficiently large $r.$ Then, by $\mf A^i$ and $\mf B^i$ we denote the extensions of, respectively, $\sf A$ and $\sf B$ to $\mc X^i_r,$ analogous to that described in Section \ref{sect31}, where, as usual, $i\in \{0,1\}$. 

    By \cite[Theorem 1.1.4]{ABHN}, $\mb u \in D(\mf A^i)$ implies $s\mapsto \|u(\cdot,s)\|_{X_x^i} \in D(\sf A).$ In particular, for $i=0$,  \begin{equation}      
    s\to \alpha_1\eta(s):=\alpha_1(s)\sup\limits_{x\in \Omega}|u(x,s)|\in L_1(\mbb R_+,dm_r).\label{etas}\end{equation} Thus, by \eqref{Assumptionofa}, 
    \begin{equation}
    D(\mf A^i) = D(\mb A^i).\label{mfA}\end{equation}
    As we shall see, under the assumptions of this section, $\mf B^i$ is a well-defined operator on $D(\mf A^i)$ for sufficiently large $r.$ 
    Then, we define the operator $\mf K^i = \mb T^i_0+\mf A^i+\mf B^i = \mf T^i+\mf B^i$ to be the restriction of the expression on the right-hand side of \eqref{TransProblem2} to $$D(\mf K^i) =D(\mb K^i) =  D(\mb T^i_0)\cap D(\mb A^i) = D(\mb T^i_0)\cap D(\mf A^i)=D(\mb T^i) = D(\mf T^i),$$
    where the `boldface' operators were defined in \eqref{dom} and $\mb K^i$ in Theorem \ref{FGeneratesSemigroup}. We do not use the results of the previous section, so we do not know \textit{a priori} that $\mb K^i$ generates a semigroup.
    
    Since we do not assume that $\beta$ satisfies \eqref{b1}, and the loss rate $\alpha_1$ is different from the gain rate $\alpha_2$, the solvability of \eqref{TransProblem2} is of independent interest, see also \cite{kerr2025c, kerr2025}. Here, however, we are mainly concerned with \eqref{TransProblem2} due to its link with \eqref{ACP1}. 
    
    Throughout this section, we always assume that (A1) is satisfied, hence, in particular, $R(\la, \mf T^i)$ is defined for $\la >0$ and $\|\la R(\la, \mf T^i)\|_{L(\mc X^i_r)}\leq 1$.    
    \subsubsection{Properties of $\beta$}
The crucial role in the proof of the generation theorem, Theorem \ref{FGeneratesSemigroup}, is played by the fact that for a fixed $c_r<1$, given by \eqref{c_r(x)Bound}, the estimate \eqref{AuEst2} allowed for keeping various constants, appearing in the calculations, below 1, and thus made possible the application of Desch's result. Unfortunately, \eqref{AuEst2} is not available for \eqref{TransProblem2} or even for \eqref{GeneralEqn} in the $\mc X_r^0$ setting. 

It turns out, however, that under a mild assumption on $\beta,$ $c_r$ decreases to 0 as $r\to \infty$, uniformly in $s$, which allows us to take arbitrarily small positive $c_r$ in \eqref{c_r(x)Bound}. Let $z = \frac{m}{s}$, $0\leq z\leq 1,$ and define the normalized moments of $\beta$ (whenever they exist) by
\begin{equation*}
\mc c_r(s) := \frac{\mc n_r(s)}{s^r}:=\frac{1}{s^r}\cl{0}{s}m^r\beta(m,s)dm= s\cl{0}{1}z^r\beta(zs,s)dz.
\end{equation*}
As for $b$, we assume that $\mc n_0(s)$ exists and there is $l\geq 0$ such that
\begin{equation}
\mc n_0(s)\leq \beta_0(1+s^l)
\label{c0}
\end{equation}
for any $s\geq 0.$ Following, e.g., \cite[Theorem 5.1.46 c)]{BLL1}, for any $s> 0$, $r\mapsto \mc c_r(s)$ is a non-increasing function, and, by the Dominated Convergence Theorem,
\begin{equation}
\lim\limits_{r\to \infty}  \mc c_r(s)=0.
\label{cry1}
\end{equation}
For our purpose, see \eqref{why}, we need this limit to be uniform in $s$, which is not always the case. 
\begin{example}\label{example1}
Consider $\beta(m,s)$ which for $s\geq 2$ is defined by
\begin{equation}
\beta(m,s) = \left\{\begin{array}{lcl}b_1(s)&\text{for}&0\leq m\leq 1,\\
b_2&\text{for}&s-1\leq m\leq s,\\
0&\text{otherwise,}&
\end{array}
\right.
\end{equation}
where $b_2<1$ is a constant and $b_1(s) = 2s(1-b_2)+b_2.$ For large particles, this model describes fragmentation in which the sizes of daughter particles are either close to the size of the parent or close to 0.  Such fragmentation processes tend to behave badly, see \cite[Example 5.1.51]{BLL1}. For such a $\beta$, 
\begin{align*}
\mc c_r(s) &= \frac{1}{s^r}b_1(s)\left(\cl{0}{1}m^rdm + b_2\cl{s-1}{s}m^rdm\right) = \frac{1}{r+1}\left(\frac{b_1(s)}{s^r} + b_2s\left(1-\left(1-\frac{1}{s}\right)^{r+1}\right)\right),
\end{align*}
and we see that \eqref{cry1} holds, but, using the l'Hospital rule,  for any fixed $r>1$,  $\lim_{s\to \infty}  \mc c_r(s)=b_2$. Hence, \eqref{c_r(x)Bound} is satisfied, but \eqref{cry1} is not uniform in $s$.  

Observe that if $b_1(s)=b_2 =1,$ then, for any  $r>1$,  $\lim_{s\to \infty}  \mc c_r(s)=1,$ and hence even \eqref{c_r(x)Bound} is not satisfied. 
\end{example}
To avoid situations described in Example \ref{example1}, we need to introduce an additional assumption. Let us recall, e.g., \cite[Theorem 4.30]{Bre} or \cite[Section 7.1]{BLL2}, that a bounded set $\mc E\in L_1(\Theta, d\mu),$ where $\mu(\Theta)<\infty$ is called uniformly integrable (or equi-integrable) if for any $\e>0$ there is $\delta>0$ such that for any measurable $\Theta_0\subset \Theta$ with $\mu(\Theta_0)<\delta$ we have
\begin{equation}
\sup\limits_{f\in \mc E}\int_{\Theta_0} |f|d\mu <\e
\label{equi}
\end{equation}
We have the following result.
\begin{proposition}\label{propei1}
Assume that there are $r_0\geq 0$ and $s_0\geq 0$ such that the set
\begin{equation}\mc E_{r_0}:= \{[0,1]\ni z\mapsto sz^{r_0}\beta(z s,s)\}_{s\geq s_0}\label{betaass}\end{equation} is equi-integrable. Then 
\begin{description}
\item a) if \eqref{aa1} and \eqref{c0} are satisfied, then $(\mf B, D(\mf A))$ is well-defined in $\mc X_r$ for any $r\geq \max\{l,r_0\},$
\item b) the limit \eqref{cry1} is uniform in $s\geq s_0$. 
\end{description}
\end{proposition}
\begin{proof}
Since $z\in [0,1]$, $\mc E_r$ is equi-integrable for any $r\geq r_0$. Then, by definition, for any $r\geq r_0$, there exists $C_r$ such that 
\begin{equation}
\sup\limits_{s\geq s_0}\frac{1}{s^{r}} \cl{0}{s} m^{r}\beta(m,s) dm =\sup\limits_{s\geq s_0}s \cl{0}{1} z^{r}\beta(zs,s) dz \leq C_r\leq C_{r_0}.
\label{equibound}
\end{equation} 
 We may assume $s_0\geq 1.$

\noindent 
a) For $\mb u\in \mc X_r^i,$ $s\mapsto \|u(s,\cdot)\|_{X_x^i} $ is measurable and we have 
\begin{equation}
\left\|\mc B\mb u\right\|_{\mc X_r^i} \leq \int_0^\infty\left( \int_m^\infty \alpha_2(s)\beta(m,s)\|u(\cdot,s)\|_{X_x^i}ds\right) w_r(m) dm.
\label{estmcB}
\end{equation}
Thus, in the calculations below, we ignore the spatial variable as it does not play any role. For a nonnegative and measurable on $\mbb R_+$ function $u$, we have, by Tonelli's theorem,  
\begin{equation}
\begin{split}
&\int_0^\infty\left( \int_m^\infty \alpha_2(s)\beta(m,s)u(s)ds\right) w_r(m) dm \\
&= 
\int_0^\infty \alpha_2(s)u(s)w_r(s) \left(\frac{1}{w_r(s)}\int_0^s w_r(m)\beta(m,s)dm\right)ds \\
&=\left(\int_0^{s_0} + \int_{s_0}^\infty\right) \alpha_2(s)u(s)w_r(s) \left(\frac{1}{w_r(s)}\int_0^s w_r(m)\beta(m,s)dm\right)ds=:I_1+I_2. 
\end{split}
\label{split}
\end{equation}
By \eqref{c0}, on $[0,s_0]$ we have 
\begin{equation}
\frac{1}{w_r(s)}\int_0^s w_r(m)\beta(m,s)dm\leq \beta_0(1+s_0^l)=: C_0.
\label{s0}
\end{equation}
Taking $r\geq \max\{r_0,l\}$, 
on $[s_0,\infty), s_0\geq 1,$ we obtain
\begin{equation}\label{beta1est}
\begin{split}
\frac{1}{w_r(s)}\int_0^s w_r(m)\beta(m,s)dm & \leq 
\beta_0\frac{w_l(s_0)}{w_r(s)} +\frac{s^{r}}{w_r(s)} \frac{1}{s^{r}} \cl{0}{s} m^{r}\beta(m,s)dm\\& \leq \beta_0+ \mc c_r(s)\leq \beta_0+C_{r_0}=C_2. 
\end{split}
\end{equation}
Using \eqref{aa1}, we have 
\begin{equation}\label{betaest}
\begin{split}
I_1+I_2 &\leq  (C_0+C_2)M\int_0^\infty \alpha_1(s)u(s)w_r(s) ds =: \bar C \int_0^\infty \alpha_1(s)u(s)w_r(s) ds,\end{split}
\end{equation}
which is finite provided $u \in D(\sf A)_+.$ Let us now return to the dependence on $x$. If $i=1$, then the thesis follows directly from \eqref{estmcB} and the comment preceding \eqref{mfA}, as the integral norm in $X_x^1$ can be moved outside the inner integral.  If $i=0,$ then, using \eqref{etas}, we observe that since $\alpha_2(s)\beta(m,s)u(x,s)\leq M\beta(m,s)\alpha_1(s)\eta(s)$ for $u\in D({\sf A})_+$, 
\eqref{betaest} yields that $(m,s)\mapsto \chi_{[m,\infty)}\alpha_2(s)\beta(m,s)\eta(s)$ is integrable on $\mbb R_+^2$ with respect to measure $dm_rds$ and hence, by Fubini--Tonelli theorem, for almost any $m$, $s\mapsto \chi_{[m,\infty)}\alpha_2(s)\beta(m,s)\eta(s)$ is integrable on $\mbb R_+$ with respect to $ds$. Therefore, by the Dominated Convergence Theorem, for $\mb u\in D(\mf A^0),$
$$
\Omega\ni x \mapsto [\mf B^0 u(x,\cdot)](m)= \int_m^\infty \alpha_2(s)\beta(m,s)u(x,s)ds \in C_0(\Omega)
$$
and hence \eqref{estmcB} and \eqref{betaest} complete the proof.

\noindent
b) In addition to \eqref{equibound}, for any $\e\in (0,1)$ we can pick up $\eta>0$ such that 
\begin{equation}
\sup\limits_{s\geq s_0}s \cl{1-\eta}{1} z^{r_0}\beta(zs,s) dz \leq \frac{\e}{2}.
\label{ei1}
\end{equation}
For these $\e$ and $\eta,$ we find $r_1>r_0$ such that $(1-\eta)^{r-r_0}C_{r_0}\leq \frac{\e}{2}$ for $r\geq r_1$. Then
\begin{equation}
\begin{split}
\mc c_r(s) &= s\cl{0}{1-\eta} z^r b(zs,s)dz +s \cl{1-\eta}{1} z^r b(zy,y) dz \\&\leq (1-\eta)^{r-r_0} s\cl{0}{1} z^{r_0}b(zs,s)dz +
s \cl{1-\eta}{1} z^{r_0}b(zs,s) dz\leq (1-\eta)^{r-r_0} C_{r_0} +\frac{\e}{2} \leq \e,
\end{split}
\label{mcc}
\end{equation}
for $r\geq r_1$ uniformly in $s\geq s_0,$ showing that $\lim_{r\to \infty}\mc c_r(s)=0$ uniformly in $s\geq s_0.$ \end{proof}
\begin{example}
Consider a homogeneous $\beta$,  
$$\beta(m,s)=\dfrac{1}{s}h\left(\dfrac{m}{s}\right),$$ which, if $\beta$ is the fragmentation kernel $b,$ describes the so called homogeneous fragmentation \cite[Section 2.2.3.2]{BLL1}. Then, assuming that $h\in L_1([0,1], z^{r_0} dz)$ for some $r_0,$
for any measurable $E\subset [0,1]$
$$s\int_{E}\beta(zs,s)z^{r_0}dz  = \int_{E} h(z)z^{r_0} dz. $$
 for $r\geq r_0$. Hence, \eqref{equi} is satisfied and the limit in \eqref{cry1} is uniform in $s$ for $s\geq s_0>0.$ 
\end{example}
\subsubsection{The generation result}
First, we provide a link between \eqref{ACP1} and \eqref{TransProblem2}. Since, as in \eqref{Tmalpha},  
$$
0\leq R(\la, \mb T^i) \leq R(\la,\mf T^i), 
$$
if for some $r\geq r_0$, $(\mf K^i, D(\mb T^i))$ generates a positive semigroup, say,  $(G_{\mf K^i}(t))_{t\geq 0}$ solving \eqref{TransProblem2} in  $\mc X^i_r$, and   
\begin{equation}
b(x,m,s)\leq \beta (m,s), \quad \text{for\;a.e.\;} x \in \Omega, (m,s)\in \Delta, 
\label{bbeta}
\end{equation}
then, by iterations,  $\mb K^i$ also generates a positive semigroup, say, $(G_{\mb K^i}(t))_{t\geq 0},$ solving \eqref{ACP1} in $\mc X^i_r,$ and satisfying  
\begin{equation}
G_{\mb K^i}(t)\leq G_{\mf K^i}(t).\label{semest}
\end{equation}
    
The main result of this section requires the following lemma.
\begin{lemma}\label{MT1}
Let \eqref{n(s)bounded}, \eqref{c0} and \eqref{betaass} be satisfied.	Then there exist $r_1> \max\{l,r_0\}$ and $\la_0$ such that for any $r\geq r_1$ and $\la>\la_0$
\begin{equation}
\|\mf B^iR(\la, \mf T^i)\|_{L(\mc X^i_r)}<1,
\label{Desc1}
\end{equation}
and hence also 
\begin{equation}
\| \mb B^iR(\la,\mb  T^i)\mb f\|_{\mc X^i_r}<1.
\label{Desc1a}
\end{equation}
\end{lemma}
\begin{proof}
The proof is similar to that of Theorem \ref{FGeneratesSemigroup} (\cite[Theorem 2.2]{BanLam20}),  but we do not assume \eqref{c_r(x)Bound}, and we cannot use \eqref{AuEst2} as we have different functions appearing in $\mc A_1$ and $\mc B_1,$ in contrast to a single $a$ in $\mc A$ and $\mc B$.  Nevertheless, we can prove \eqref{Desc1}. 
Denote
$$v(s,\mb f):=\|[R(\lambda, \mf T^i)\mb f](\cdot,s)\|_{X_x^i}.$$
Applying \eqref{AR} to the current setting, for $\mb f\in \mc X^i_r$ we obtain   
				\begin{equation}
                ||R(\lambda,\mf T^i)\mb f||_{\mc {X}^i_r} = \int\limits_0^\infty v(s,\mb f)ds_r\leq  \frac{1}{\lambda}||\mb f||_{\mc {X}^i_r},
\label{Rb1}
\end{equation}
and, using the fact that the second and third terms of \eqref{alph1est} lead to \eqref{AuEst2} for $\mb A$ independent of $x$, as is the case with $\mf A^i$, we get 
\begin{equation}
\label{Rb2}
\begin{split}
	\|\mf A^iR(\lambda,\mf T^i)\mb f||_{\mc {X}^i_r}=\int\limits_0^\infty \alpha_1(s)v(s,\mb f) ds_r 	&\leq  ||\mb f||_{\mc {X}^i_r}.\end{split}\end{equation}			Passing to the main part of the proof, let $\mb f\in \mathcal{X}^i_{r,+}$ and $\lambda>0$. Using \eqref{estmcB}, we get
\begin{equation*}
\|\mf B^iR(\lambda,\mf T^i)\mb f||_{\mc X^i_r}  \leq  \int_0^\infty\left( \int_m^\infty \alpha_2(s)\beta(m,s)v(s,\mb f)ds\right) w_r(m) dm =I_1+I_2,
									\end{equation*}
where the split is as in \eqref{split}.  Setting $\alpha_0 =\esssup\limits_{s\in[0,s_0]}\alpha_2(s)$, from \eqref{s0} we obtain 	\begin{equation*}
I_1\leq \frac{\alpha_0 C_0}{\la}\|\mb f\|_{\mc X^i_r}.
\end{equation*}
To obtain a suitable estimate for the integral over $[s_0,\infty)$, we use Proposition \ref{propei1} to chose $r_1>\max\{r_0,l\}$ such that for $r\geq r_1,$  $\mc c_r(s)$ is  small enough for 
\begin{equation}
\beta_0 \mc w_r(s)+\mc c_r(s)< \frac{1}{2M}, \quad s \in [s_0,\infty),\label{why}
\end{equation}
where $M$ was defined in \eqref{aa1}. Then, using \eqref{Rb2},
\begin{equation*}
\begin{split}
I_2& < \frac{1}{2M}\int\limits_{s_0}^{\infty}\!\alpha_2(s)v(s,\mb f)ds_r<\frac{1}{2}||\mb f||_{\mc X^i_r}.
					\end{split}
				\end{equation*}
				Now, choosing $\la$ large enough for $\frac{\alpha_0 C_0}{\la}<\frac{1}{2}$, we get
$$||\mf B^iR(\lambda,\mf T^i)\mb f||_{\mc X^i_r} =I_1+I_2<\left(\frac{\alpha_0 C_0}{\lambda}+\frac{1}{2}\right)\|\mb f\|_{\mc X^i_r}<\|\mb f\|_{\mc X^i_r}.$$ 
				\end{proof}
This immediately gives the generation in $\mc X^1_r$ due to the Desch theorem, as in Theorem \ref{FGeneratesSemigroup}. 
\begin{corollary}\label{corg1}
 Assume that the assumptions of Lemma \ref{MT1} are satisfied. Then $(\mf K^1,D(\mb T^1))$ generates a positive $C_0$-semigroup, say $(G_\mf K^1(t))_{t\geq0}$, on $\mathcal{X}^1_r$. Hence also $(\mb K^1,D(\mb T^1))$ generates a positive $C_0$-semigroup on $\mc X^1_r.$
 \end{corollary}
 Another case, when Lemma \ref{MT1} gives the generation is if $(G_{\mb T}(t))_{t\geq 0}$ is an analytic semigroup. Then, the Arendt--Rhandi theorem, \cite[Theorem 1.1]{AR} or \cite[Theorem 4.1.1]{BLL1} (where the last statement follows as the assumptions in \textit{op. cit.} can be applied directly to \eqref{Desc1a}), yield  
 \begin{corollary}\label{corg1a}
 Assume that the assumptions of Lemma \ref{MT1} are satisfied and $(G_{\mb T}(t))_{t\geq 0}$ is an analytic semigroup. Then $(\mf K,D(\mb T))$ generates a positive analytic $C_0$-semigroup, say $(G_\mf K(t))_{t\geq0}$, on $\mathcal{X}_r$. Hence also $(\mb K,D(\mb T))$ generates a positive analytic $C_0$-semigroup on $\mc X_r.$
 \end{corollary}
Unfortunately, in general,  \eqref{Desc1} is insufficient for the application of the Desch theorem.  Hence, we need a stronger assumption. First, however, we shall prove an auxiliary result which is also of independent interest. 
\begin{proposition} Let the assumptions of Lemma 3.2 be satisfied and $r_1$ be as defined there.   
For any $r\geq r_1,$ the operators $({\sf F},D({\sf  A}))= ({\sf A}+{\sf B}, D({\sf A}))$ and $(\mf F, D(\mf A)) = (\mf A +\mf B, D(\mf A))$ generate positive analytic semigroups, say, $(G_{{\sf F}}(t))_{t\geq 0}$ and $(G_{\mf F}(t))_{t\geq 0}$ in, respectively, $X_{m,r},$ for $m\in \mbb R_+$,  and $\mc X_r,$ satisfying 
\begin{equation}
[G_{\sf F}(t) f(x,\cdot)](m) = [G_{\mf F}(t) \mb f](x,m), \quad \mb f \in \mc X_r,\, \quad x\in \Omega.
\label{frageq}
\end{equation}
\end{proposition}
\begin{proof}
The generation of an analytic semigroup by $\sf F$ is an immediate consequence of Lemma \ref{MT1} applied to the case with $\mc T_0=0$ in $X_{m,r},$ and of the Arendt--Rhandi theorem. The extension to $\mc X^1_r$ follows as in Corollary \ref{corg1}. To extend the result to $\mc X_r^0$, first we observe that since 
$$[{\sf B}R(\la,{\sf A})u(x,\cdot)](m) = \int_m^\infty \frac{\alpha_2(s)}{\la+\alpha_1(s)}\beta(m,s) u(x,s)ds,$$
$({\sf B}R(\la,{\sf A}))^n$ is an integral operator with $x$-independent kernel, and, by,  e.g., \cite[Theorem 5.10]{BaAr},
\begin{align*}
\|R(\la, \mf F^0)\mb u\|_{\mc X_r^0} &\!= \!\int^\infty_0\!\sup\limits_{x\in \Omega}|[R(\la,{\sf F}) u(x,\cdot)](m)|dm_r \!=\! \int^\infty_0 \!\left\|R(\la,{\sf A})\sum\limits_{n=1}^\infty ({\sf B}R(\la,{\sf A}))^n\mb u\right\|_{X_x^0}\!\! (m)dm_r\\
&\leq \int^\infty_0 \frac{1}{\la+\alpha_1(m)}\sum\limits_{n=1}^\infty \left(\left[({\sf B}R(\la,{\sf A}))^n\|\mb u\|_{X_x^0}\right] (m)dm_r\right)\\
&= \|R(\la, {\sf F})\|\mb u\|_{X_x^0}\|_{X_{m,r}}\leq \|R(\la, {\sf F})\|_{L(X_{m,r})} \|\mb u\|_{\mc X_r^0},
\end{align*}
hence $R(\la, \mf F^0)$ satisfies the Hille--Yosida (and sectorial) estimates as $R(\la, {\sf F})$ does. Since $D(\mf A^0)\subset D(\mf A^0) = D(\mb A^0)$, the density of the domain is obvious, and $(\mf F^0,D(\mb A^0))$ generates an analytic semigroup. Equality \eqref{frageq} follows from obvious equality
$$
[\mf B R(\la, \mf A^i)\mb f](x,m) = [{\sf B} R(\la, {\sf  A}) f(x,\cdot)](m), \quad \mb f \in \mc X_r,\, \text{for\;a.e.}\;x\in \Omega,
$$ and the above representation of $R(\la, \sf F)$.
\end{proof}
The main purpose of the next corollary is to prove the existence of the advection--fragmentation semigroup in $\mc X_r^0,$ but, as later we shall need the commutativity property \eqref{comm1} in both cases, we formulate the result in general setting. 
\begin{corollary}\label{corg2}
Let all assumptions of Lemma \ref{MT1} be satisfied and, in addition, let $\mc T_0$ be independent of $m$. Then $(\mf K,D(\mb T))$ generates a positive $C_0$-semigroup, say $(G_\mf K(t))_{t\geq0}$, on $\mathcal{X}_r$. Hence also $(\mb K,D(\mb T))$ generates a positive $C_0$-semigroup on $\mc X_r.$
\end{corollary}
\begin{proof}
Let $\mc X_r\ni \mb f = \phi f, \phi \in X_x, f\in L_1(\mbb R_+,dm_r)$. Then
\begin{align*}
        [ G_{\mb T_0}(t)( G_{\mf F}(t)\mb f)](x,m) &= [G_{\mb T_0}(t) \mb f](x)[ G_{\mf F}(t)]\phi (m) = [G_{\mf F}(t)(G_{\mb T_0}(t)\mb f)](x,m),\quad t\geq 0.
        \end{align*}
Now, as in Proposition \ref{glueprop1}, we use the fact that $\mc X_r =  L_1(\mbb R_+,dm_r)\hat \otimes_\pi X_x$ and the linearity, to claim that the equality is valid for any $\mb f\in \mc X_r$. Next, by \cite[Sections I 5.15 \& II 2.7]{EN}, $(G(t))_{t\geq 0} := (G_{\mb T_0}(t) G_{\mf F}(t))_{t\geq 0}$ is a $C_0$-semigroup in $\mc X_r,$ whose generator restricted to $D(\mb T_0)\cap D(\mf F)$ is $\mb T_0+{\mf F}$ and $D(\mb T_0)\cap D(\mf F) = D(\mb T_0)\cap D(\mb A)$ is its core.  Since, however, $(\mf K,D(\mb T)) = (\mb T_0+\mf F, D(\mb T_0)\cap D(\mb A)$ has the resolvent $R(\la, \mf K)$ for large $\la$, it must be the generator and 
\begin{equation}
G_{\mf K}(t)\mb f =  G_{\mb T_0}(t) G_{\mf F}(t)\mb f =  G_{\mf F}(t) G_{\mb T_0}(t)\mb f, \quad \mb f \in \mc X_r.
\label{comm1}
\end{equation}
\end{proof}
                \subsubsection{Application of the analyticity of the fragmentation operator}
If  $\mc T_0$ is a diffusion operator, then $\sem{G_{\mb K}},$ in the cases described by Corollary \ref{corg1} or \ref{corg2}, can be proved to be an analytic semigroup by the Arendt--Rhandi theorem, \cite[Theorem 1.1]{AR}, mentioned above. However, the identification of interpolation spaces between $D(\mb K)$ and $\mc X_r,$ needed for moment regularisation, is far from obvious. We can use, however, the analyticity of the fragmentation semigroup to prove necessary estimates if $\mc T_0$ is independent of $m$.

Since $({\mf F}, D(\mf A))$ generates an analytic semigroup in $\mc X_r$ for $r\geq r_1$ for some $r_1>1$, we can apply the theory developed in \cite{BanDCDS20, BaLa12a}, see also \cite{kerr2025}. Referring the reader to \textit{op.cit.} for details, here we mention that, defining 
$$
\mf F_{\omega}:= \mf F - \omega I  = \mf A_{\omega} + \mf B = 
 \mf A - \omega I+\mf B,
$$
where $\omega>1$ is greater than the type of $(G_{\mf F}(t))_{t\geq 0}$ and using the fact that $D(\mf F_{\omega}) = D(\mf A_{\omega})$, we identify the real interpolation space $D_{\mf F_{\omega}}(\mu, 1)$, see \cite[Corollary 2.2.3]{Lun}, for $\mu \in [0,1]$  with
\begin{equation}
 \mc X_{r}^{(\mu)} := \left\{ \mb f \in \mc X_{r}:\; \ \int_0^\infty \|\mb f(m)\|_{X_x}(\omega + \alpha_1(m))^\mu \,d m_r < \infty \right\}.
 \label{frps1}
\end{equation}
 In general, by \cite[Proposition 2.2.9]{Lun} (see the proof with the correct range of parameters), for an analytic semigroup  $\sem{G_{\mf F}}$  on $\mc X_r$, there are constants $\omega_r, M_r^{\mu}$ such that  for $\mu \in [0,1],$
\begin{equation}\label{eq2.4a}
\|G_{\mf F}(t)\mb f\|_{\mc X_r^{(\mu)}}\leq \frac{M_{r}^{(
\mu)}e^{\omega_{r}t}}{t^\mu}\|\mb f\|_{\mc X_r}, \quad t>0.
\end{equation}
 
Explicitly, \eqref{eq2.4a} expresses a moment improving property of $\sem {G}$  at the cost, however, of worsening the regularity at $t=0.$ 
\begin{theorem}\label{th34}
Let all assumptions of Lemma \ref{MT1} be satisfied and, in addition, let $\mc T_0$ be independent of $m$.  Then, for any $r\geq r_1$
and $\mu\in [0,1]$, there are constants $M_r^{(\mu)}$ and (independent of $\mu$) $\omega_r$ such that for any $t>0, \mb f\in \mc X_r,$
\begin{equation*}
\begin{split}
\|G_{\mb K}(t)\mb f\|_{\mc X^{(\mu)}_r}&\leq  \|G_{\mf K}(t)\mb f\|_{\mc X^{(\mu)}_r} \leq  \frac{M_r^{(\mu)}}{t^\mu}e^{\omega_r t}\|\mb f\|_{\mc X_r}.
\end{split}
\end{equation*}
If \eqref{aCondition} is satisfied, then for $q:=\mu\gamma \leq \gamma$
\begin{equation}
    \| G_{\mb K}(t) \mb f \|_{\mc X_{r+q }} \leq \| G_{\mf K}(t) \mb f \|_{\mc X_{r+q}}\leq C_1\|G_{\mf K}(t) \mb f \|_{\mc X_r^{(\mu)}} \leq C_1\frac{M_r^{(\mu)}}{t^\mu}e^{\omega_r t}\|\mb f\|_{\mc X_r}, \quad \mb f \in \mc X_r.
    \label{momimp}
\end{equation}
\end{theorem}
\begin{proof}
Since $\sem{G_{\mf F}}$ is analytic, \eqref{semest}, \eqref{comm1}, and \eqref{eq2.4a}, yield for $\mb f\in \mc X_r$ and $\mu \in [0,1]$, 
\begin{equation*}
\begin{split}
\| G_{\mb K}(t)\mb f\|_{\mc X^{(\mu)}_r}&\leq  \| G_{\mf K}(t)\mb f\|_{\mc X^{(\mu)}_r} = \| G_{\mf F}(t)G_{T_0}(t)\mb f\|_{\mc X^{(\mu)}_r} \leq \frac{M_r^{(\mu)}}{t^\mu}e^{\omega_r t}\| G_{\mb T_0}(t)\mb f\|_{\mc X_r}\\
&= \frac{M_r^{(\mu)}}{t^\mu}e^{\omega_r t}\|G_{\mb T_0}(t)\mb f\|_{\mc X_r}\leq \frac{M_r^{(\mu)}}{t^\mu}e^{\omega_r t}\|\mb f\|_{\mc X_r}, 
\end{split}
\end{equation*}
where we used the fact that $\sem{G_{\mb T_0}}$ is contractive. 

For the last statement, using the obvious estimates $$
 \frac{1+m^{r+q}}{(1+m^q)(1+m^r)}\leq 1, \quad m\geq 0,
\quad 
\frac{1+m^q}{(\omega + \alpha_1(m))^{\frac{q}{\gamma}}}\leq \frac{1+m_0^q}{\omega^{\frac{q}{\gamma}}} =: C_{m_0}, \quad m\in [0,m_0],
$$
we get  
\begin{align*}
&\int\limits_0^\infty \|\mb f(m)\|_{X_x}\,d m_q \leq 
\int\limits_0^\infty \|\mb f(m)\|_{X_x}(1+m^q)\,d m_r \\& \leq C_{m_0}\int\limits_0^{m_0} \|\mb f(m)\|_{X_x}(\omega+\alpha_1(m))^{\frac{q}{\gamma}}\,d m_r + \max\{\omega^{-\frac{q}{\gamma}}, a_0^{\frac{q}{\gamma}}\}\int\limits_{m_0}^\infty \|\mb f(m)\|_{X_x}\left(\omega^{\frac{q}{\gamma}}+\alpha^{\frac{q}{\gamma}}_1(m)\right)\,d m_r\\
& \leq C_1\int\limits_{0}^\infty \|\mb f(m)\|_{X_x}(\omega+\alpha_1(m))^{\frac{q}{\gamma}}\,d m_r = C_1\|\mb f\|_{\mc X_r^{\left(\frac{q}{\gamma}\right)}},
\end{align*}
where $C_1 = \max\{C_1,\omega^{-\frac{q}{\gamma}}, a_0^{\frac{q}{\gamma}}\}$. Since $\mu = \frac{q}{\gamma},$ 
\begin{align*}
    \|G_{\mf K}(t) \mb f \|_{\mc X_{r+q}} \leq C_1\|G_{\mb K}(t) \mb f \|_{\mc X_r^{(\mu)}} \leq C_1\frac{M_r^{(\mu)}}{t^\mu}e^{\omega_r t}\|\mb f\|_{\mc X_r}, \quad \mb f \in \mc X_r,
\end{align*}
from which \eqref{momimp} follows immediately  by \eqref{semest}.
\end{proof}
\section{Applications}\label{sec4}
\subsection{Advection--fragmentation equation}\label{sec41}
In this section, $\mc T_0$ is the advection operator describing the movement of the particles in the physical space $\mbb R^N$ due to the flow of the surrounding medium. To avoid being bogged down by technicalities, we consider the flow on $\mbb {R}^N$. The case of bounded domains requires delicate handling of the boundary conditions, see \cite[Section 10.3]{BaAr} or \cite{ABL3}, but, to a large extent, can be dealt with by the approach presented here as long as the semigroup is explicitly given by the composition of the initial state with the flow.

We assume that the velocity field  $\mb \omega(x,m)$ of the moving medium satisfies   
	\begin{itemize}
	\item[(a1)] $\mb \omega:\mbb R^N\times\mathbb{R}_+\rightarrow\mbb R^N$ is independent of time, globally Lipschitz continuous with respect to $x$ uniformly in $m$, and
    \subitem (i) uniformly continuous with respect to $m$ uniformly in $x$ if $i=0$, or 
    \subitem (ii) globally Lipschitz continuous with respect to $m$ uniformly in $x$ if $i=1$. 
    
    We can assume that the Lipschitz constant is the same and denote it by $\kappa>0$.
	\item[(a2)] for any $m\in \mathbb{R}_+$,  $\mb \omega$ is divergence--free with respect to $x$.
\end{itemize} 
Thus, for $u \in C_c^1(\mbb R^N)= \{u\in C^1(\mbb R^N):\;u\;\text{is\;compactly\;supported}\}$, $\mc T_0$ can be defined by \begin{equation}
				[\mathcal{T}_0u](x,m)= -\mb \omega(x,m)\cdot \nabla_x u(x,m),\quad (x,m)\in \mbb R^N\times \mbb R_+.\label{ops1}
                \end{equation}
 Following Section \ref{sect31}, for  $i\in\{0,1\}$ and  $m \in \mbb R_+,$ we define the advection operator $(T^i_{0,m}, D(T^i_{0,m}))$ as the extension of $(\mathcal{T}_{0}, C^1_c(\mbb R^N))$ that generates a $C_0$-semigroup in $X_x^i$, $i=1,2$, as defined in the theorems below. 
	
    Let us define the flow  $\mb\phi$ of the system 
	\begin{equation}\label{IVP1}
		\begin{cases}
			\frac{d}{ds} y(s)&=\mb \omega( y(s),m), \quad s\in\mathbb{R},\\
			y(t)&=x,
		\end{cases}
	\end{equation}
    so that $\mb y(s) = \mb\phi(x,s-t, m) ,$ \cite[Proposition 10.1]{BaAr}. 	The global existence and uniqueness of solutions to \eqref{IVP1} follow from assumption (a1). Next,  as in \eqref{TransProblem}, we define 
    \begin{equation}\label{TransProb2*}
		\begin{cases}
			T^i_{m}u&=\mathcal{T}u \\
			D(T^i_{m})&=  D(T^i_{0,m}),
		\end{cases}
	\end{equation}
    where $a$ satisfies the relevant version of assumption \eqref{Assumptionofa}. 
           Here, the semigroups $\sem{G_{T^i_m}}$ are explicitly known, thus we can skip some of the technicalities of Section \ref{sect31}. 
           Indeed, we have,
	\begin{theorem}\cite[Theorem 10.4]{BaAr}\label{TransportGeneration}
		For  $m\in \mbb R_+$, the family $(G_{T^1_{m}}(t))_{t\geq 0}$ defined by
		\begin{equation}\label{ResolventofDtilde}
			[G_{T^1_{m}}(t)f](x,m)=e^{-\int_0^t a(\mb\phi(x,s,m))ds}f(\mb\phi(x,t,m))
		\end{equation}  for any $f\in {X}^1_x$ and $t\geq 0$, is a substochastic semigroup, generated $T_{m}^1 = \mc T$ on the domain $D(T_{m,0}^1) = \{u\in X_{x}^1:\;\mc T_0u \in X_{x}^1\},$ where $\mc T_0u$ is defined in the sense of distributions. 
	\end{theorem}
    \begin{remark}
    Using the argument from \cite[Section II.3.28]{EN}, one can prove that the set of compactly supported Lipschitz functions is the core of $T^1_{0,m}$ and thus $T^1_{0,m}$ can also be characterized as the closure of $\mc T_0$ restricted to this set.   
    \end{remark}
	Recalling the definitions \eqref{dom}, we have  
           \begin{corollary}\label{corsubs}
Let us fix $r\geq 0$. The family of operators, defined for any $\mb f \in \mc X^1_r$ by
\begin{equation}\label{GXr}
			[G_{\mb T^1}(t)\mb f](x,m)=e^{-\int_0^t a(\mb\phi(x,s,m))ds}\mb f(\mb\phi(x,t,m),m), \quad t\geq 0, x \in \mbb R^N, m\in \mbb R_+,
		\end{equation} 
        is a substochastic semigroup on $\mc X^1_r$ generated by $\mb T^1=\mc T$ on $$D(\mb T^1) = D(\mb T^1_0)\cap D(\mb A^1).$$
        If $a\equiv 0,$ then $\sem{G_{\mb T^1}}$ is stochastic. 
    \end{corollary}
    \begin{proof}
By Proposition \ref{prop31} and Theorem \ref{TransportGeneration}, it suffices to prove the statement in the case $a\equiv 0$. Since the semigroups $\sem{G_{T^1_{0,m}}}$ are known, we can use Proposition \ref{ParameterSemigroup}. Thus, we need to show that for each $t\geq 0$ and $\mb f=f(\cdot,\cdot)\in \mc X^1_r,$
  \begin{equation}
  \mbb R_+\ni m\mapsto \Psi(t,\cdot,m)=f(\mb\phi(\cdot,t,m),m)\in X_x
  \label{Phidef}
  \end{equation}
  is a Bochner-measurable $ X^1_x$-valued function. First, we observe that $\Psi$ is a scalar measurable function on $\mbb R^N\times \mbb R_+$, since the inner function in $\Psi,$ $(\mb z, m) = \mb \Phi(x, m):= (\mb\phi(x,t,m),m)$ is bi-Lipschitz. This follows from  
  $\mb \Phi^{-1}(\mb z,m) = (\mb\phi(\mb z,-t,m),m),$ and the assumption (a1)(ii) that gives, by the Gr\"{o}nwall inequality, 
  \begin{equation}
  \|\mb \phi (x,t,m_1)  - \mb \phi (y,t,m_2)\|\leq (\|x-y\|+t|m_1-m_2|)e^{\kappa t}.
  \label{GI1}
  \end{equation}
    Now, by assumption (a2) and \eqref{Phidef}, we have for $\mb f \in \mc X^1_r$, 
  \begin{equation}\label{divfree}
  \begin{split}
           \int\limits_{0}^{\infty}  \int\limits_{\mbb R^N}\left|\Psi(t,\cdot,m)\right| dx dm_r &\leq 
      \int\limits_{0}^{\infty}  \left(\,\int\limits_{\mbb R^N} |f(\mb z,m)|d\mb z\right) dm_r = \|\mb f\|_{\mc X^1_r}<\infty,
      \end{split}
  \end{equation}
  and the identification \eqref{bb} yields  
  $
   \mbb R_+\ni m\mapsto \Psi(t,\cdot,m) \in L_1(\mbb R_+, X^1_x, dm_r).
  $
  Hence, in particular, it is a Bochner-measurable $X^1_x$-valued function. 

Thus, assumption (A1) is satisfied and hence, by Proposition \ref{prop31}, the semigroup $\sem{G_{\mb T^1}}$ is a substochastic semigroup generated by $\mb T^1$  on $D(\mb T^1)=D(\mb T^1_0)\cap D(\mb A^1)$. 
  
   To prove the last statement, we note that if $a\equiv 0,$ then the assumption that $\mb \omega$ is divergence-free implies, by direct integration, see \cite[Theorem 10.2]{BaAr}, that for almost every $m$ 
   $$
   \|G_{T^1_{0,m}}(t)\mb f(\cdot,m)\|_{X^1_x} = \|\mb f(\cdot,m)\|_{X^1_x}, \quad \mb f \in \mc X^1_r, 
   $$
   and the statement follows by integration of the above formula with respect to $dm_r$.  
    \end{proof}
An analogous result in $X^0_x=C_0(\mbb R^N)$ seems to be folklore; see \cite[Section II.3.28]{EN}. We will fill in some details below and extend it to $\mc X_r^0$.      \begin{theorem}Let us fix $r\geq 0$.\begin{enumerate}
\item For any $m\in \mbb R_+$, the formula 
\begin{equation}
[ G_{T^0_{m}}(t)f](x,m)=e^{-\int_0^t a(\mb\phi(x,s,m))ds}f(\mb \phi(x,t,m)),\quad t\geq 0,\,(x,m) \in \mbb R^N\times \mbb R_+,
\label{C0sem}
\end{equation}
for $f \in C_0(\mbb R^N)$, defines a substochastic $C_0$-semigroup in $X_x^0.$ If $\mb \omega$ is continuously differentiable with bounded derivative on $\mbb R^N$, then the  generator $T^0_{m}$ is the closure of $\mc T_{m}$ defined on $C_c^\infty(\mbb R^N).$ 
\item  The family of operators, defined for any $\mb f \in \mc X^0_r$ by
\begin{equation}\label{GXr0}
			[G_{\mb T^0}(t)\mb f](x,m)=e^{-\int_0^t a(\mb\phi(x,s,m))ds}\mb f(\mb\phi(x,t,m),m), \quad t\geq 0,\, (x,m) \in \mbb R^N\times \mbb R_+,
		\end{equation} 
        is a substochastic semigroup on $\mc X^0_r$ generated by $\mb T^0=\mc T$ on $$D(\mb T^0) = D(\mb T^0_0)\cap D(\mb A^0).$$
        \end{enumerate}
\end{theorem} 
\begin{proof}
 1. We can focus on the proof that   $[G_{T^0_{0,m}}(t)f(\cdot,m)](x) := f(\mb \phi(x,t,m))$ is a $C_0$-semigroup in $C_0(\mbb R^N)$ for any $m\in \mbb R_+$. For, once it is established, then, by \eqref{aa0c}, $u\mapsto a(\cdot,m) u$ is bounded on $X_x^0$ for any $m\in \mbb R_+$ and hence the existence of $\sem{G_{T^0_{m}}}$ follows from the Bounded Perturbation Theorem and \eqref{C0sem} follows for smooth $f \in D(T^0_{m})$ by the unique solvability of the initial value problem $\p_t u = \mc T u, \quad u(0) = f,$ and passing to the limit in $X_x^0$, which is not affected by the multiplication by $e^{-\int_0^t a(\mb\phi(x,s,m))ds}.$ 

Due to assumption (a1), $x\mapsto \mb\phi(x,t,m)$ is a globally bi-Lipschitz mapping from $\mbb R^N$ onto $\mbb R^N$ for any fixed $t,m$. Hence, 
$$\lim\limits_{\|x\|\to \infty} f(\mb \phi(x,t,m),m) = \lim\limits_{\|y\|\to \infty} f(y,m)= 0,
$$
and the mapping $G_{T^0_{0,m}}(t)$ leaves $C_0(\mbb R^N)$ invariant for any $t\geq 0$. Let us skip the dependence on $m$ for a moment, as it is not relevant. First, the forward and backward Gr\"{o}nwall inequalities give, for any $t\in \mbb R,$
\begin{equation}\|\mb \phi(t,x)-x\|\leq |t|\|\mb\omega(x)\|e^{L|t|}, \quad x\in \mbb R^N.\label{unicont}
\end{equation}
 Since the formula for $G_{T^0_{0,m}}(t)$ is the same as in $X^1_x$, the semigroup property follows. To show the strong continuity, we observe first that any $f \in C_0(\mbb R^N)$ is uniformly continuous. This follows as $f$ is uniformly continuous on any closed ball and uniformly small outside it. 
 Let $f\in C_0(\mbb R^N)$. For any $\e,$ there is $R$ such that  $|f(x)|\leq \frac{\e}{2}$ for $\|x\|\geq R.$  Then, with $\omega_R: = \sup_{x\in \ov B(0,R)}\|\mb \omega(x)\|$, by \eqref{unicont}, 
 \begin{equation}
 \mb \phi(t,\ov B(0,R))\subset \ov B\left(0, R+ te^{Lt}\omega_R\right)\subset \ov B\left(0, R+ e^{L}\omega_R\right) =:K_R.\label{KR}
 \end{equation}
Since $K_R$ is compact, for any $\e>0$ there is $\delta$ such that $|f(x)-f(y)| <\e,$ whenever $x,y \in K_R$ and $\|x-y\|\leq \delta.$ Since $t\mapsto te^{Lt}=:\eta(t)$ is a continuous increasing function with $\eta(0)=0,$ there is $t_\delta<1$ such that for $t\in [0,t_\delta]$ we have $\|\mb\phi(x,t)-x\|\leq \delta$ for $x\in K_R.$

Now, we have the following possibilities. First, if both $x$ and $\mb\phi(x,t)$ are in $K_R$ for some $t\in [0,t_\delta],$
then $|f(\mb \phi(x,t)) - f(x)|<\e$ on account of the uniform continuity of $f$ on $K_R$. If $x\in K_R$ and  $\mb\phi(x,t)\notin K_R$ for some $t\in [0,t_\delta]$, then $x\notin \ov B(0,R)$ by \eqref{KR},
as the flow cannot reach the outside of $K_R$ from $\ov B(0,R)$ in time less than 1. Thus, $x\notin \ov B(0,R)$ and  $\mb\phi(x,t)\notin \ov B(0,R)$ and $|f(\mb \phi(x,t)) - f(x)| \leq \e$. Therefore, if $x\notin K_R$, then $\mb\phi(x,t)\notin \ov B(0,R)$ for $t \in [0,t_\delta]$ as in such a case, $\ov B(0,R)\ni y = \mb\phi(x,t')$ for some $t'\in (0,t_\delta],$ hence  $x = \phi(y,-t')$ for some $y\in \ov B(0,R)$, which is impossible by \eqref{unicont} and \eqref{KR}, and, as above, $|f(\mb \phi(x,t)) - f(x)| \leq \e$. Finally, if both $x$ and $\mb\phi(x,t)$ are outside $K_R$, then obviously $|f(\mb \phi(x,t)) - f(x)| \leq \e.$

This shows that $\sem {G_{T^0_{0,m}}}$ is a strongly continuous semigroup on $X^0_x$. The proof of the characterisation of the generator can be found in \cite[Section II.3.28]{EN}.

\noindent
2. Now, we revert to the dependence of the involved functions on $m$ to show that $m\to [ G_{T^0_{0,m}}(t)f](\cdot,m)$ is measurable as an $X_x^0$ valued function and hence the corresponding 
$\sem{G_{\mb T^0_0}}$ is a stochastic semigroup on $\mc X^0_r = L_1(\mbb R_+, X_x^0,dm_r)$. 

 Since $\mb f\in\mc X_r^0$, there are characteristic functions $\chi_{i,n}(m)$ and $X^0_x$-functions $ f_{i,n}$ such that 
$$
\mb f(m)= \lim\limits_{n\to \infty} \mb f_n(m)= \lim\limits_{n\to \infty} \sum\limits_{i=1}^n\chi_{i,n}(m)f_{i,n} 
$$ 
in $X^0_x$ for almost any $m$. Now, 
$$
[G_{T_{0,m}}(t)f_n(\cdot,m)](x) =  \sum\limits_{i=1}^{n}\chi_{i,n}(m)f_{i,n}(\mb \phi(x,t,m))
$$
and
$$
\sup\limits_{x\in \mbb R^N} \left|f(\mb\phi(x,t,m),m)-\sum\limits_{i=1}^n\chi_{i,n}(m)f_{i,n}(\mb \phi(x,t,m))\right| \leq 
\sup\limits_{z\in \mbb R^N} \left|f(z,m)-\sum\limits_{i=1}^{n}\chi_{i,n}(m)f_{i,n}(z)\right|.
$$
To conclude the proof, it suffices to show that for any fixed $t\geq 0,$ $m\mapsto f(\phi(x,t,m))$ is measurable for any $f\in X^0_x.$ In fact, we prove that it is a continuous $X^0_x$-valued function. 
We mentioned above that $f$ is uniformly continuous on $\mbb R^N$. 
Using \eqref{GI1} with $x=y$, we see that for any $\delta>0$ and any $t\geq 0,$ there is $\eta>0$ such that if $|m_1-m_2|<\eta,$
 $$
  \sup\limits_{x\in \mbb R^N}\|\mb \phi (x,t,m_1)  - \mb \phi (x,t,m_2)\|\leq \delta,
  $$
 hence 
 $$\sup\limits_{x\in \mbb R^N}|f(\mb\phi(x,t,m_1))-f(\mb\phi(x,t,m_2))|<\e$$
 and $m\mapsto f(\mb \phi(\cdot,t,m))$ is $X_x^0$-continuous. Therefore, $m\mapsto \mb f(m)$ is measurable. Hence, $\sem{G_{\mb T^0_0}}$ is a $C_0$-semigroup on $\mc X^0_r$ and statement 2. follows by Proposition \ref{prop31}. \end{proof}

   \subsection{Diffusion--fragmentation} \label{sec42}
    
The theory of the diffusion equation in $L_1(\Omega,dx)$ and $C(\ov{\Omega})$ spaces is complicated due to problems with identifying the domain of the generator and thus proving that it is independent of the coefficients. We consider a case when it is possible. The results in this section are based on ideas from  \cite[Section 4]{BanTTSP2007}.

Let us consider the diffusion equation, 
\begin{equation}\begin{split}
\p_t u(t,x,m) &=[\mc{T_0}u](t,x,m) =\nabla_x\cdot(d(x,m)\nabla_x u(t,x,m)),\\
u(0,x,m)&=\mr{u}(x,m),\end{split}\label{difeqd}
\end{equation}
in $\Omega \times \mbb{R}_+$, where $\Omega\subseteq \mbb R^N$ is a bounded open set with a $C^2$  boundary, and for 
almost any $m\in \mbb R_+$,  $d(\cdot,m)\in C^1(\overline{\Omega})$ and there exist $d_{\min}(m)>0$ and $d_{\max}(m)<\infty$ such that 
\begin{equation}
 d_{\min}(m)\leq
d(x,m)\leq d_{\max}(m), \quad x\in \Omega.
\label{dass1}
\end{equation}
A natural boundary condition for \eqref{difeqd} is 
\begin{equation}
\left.\partial_{\mb n} u\right|_{\p\Omega}= 0, \label{nbc}
\end{equation}
where $\p_{\mb n}$ is the outward normal derivative at the boundary. 

In this part, we do not use the dependence on $m,$ and thus we drop it from the
notation, remembering, however, that all constants can depend on $m$.  The
transversal derivative associated with $\mc {T}_0$ coincides with the normal derivative of \eqref{difeqd}, and hence we can use \cite{Am} or \cite{Fa}. 
 Let $i\in\{0,1\}.$ First, we denote by $T^i_0$  the closure of the restriction of $\mc {T}_0$ to $\{u \in
C^2(\ov{\Omega});\;\p_{\mb n}u =0\;{\rm on}\;\p\Omega\}$ in, $X_x^i$,  see \cite[Lemma 9.1]{Am}. Then, \cite[Theorems 8.2 \& 10.3]{Am} or \cite[Theorem 4.8.3]{Fa},  $(T^i_0, D(T_0))$, generate compact, analytic and stochastic semigroup on respective $X^i_x$. Moreover, we have $D(T^1_0) \subset W^{1}_q(\Omega)$ for all $q\in \left [1,\frac{3}{2}\right)$, \cite[Proposition 9.2]{Am}, and $D(T^0_0) \subset \bigcap_{q\geq 1}W^{1}_q(\Omega)\subset C^1(\ov{\Omega})$, by the Sobolev embedding theorem.   Thanks to this, we can prove that $D(T^i_0) = D(\Delta^i),$ where $\Delta^i$ is the realisation of the Laplacian which generates a semigroup on the respective $X^i_x$.

Indeed, let $u_n\to u$ and $\mc T_0u_n \to f$ as $n\to \infty$ for some sequence $(u_n)_{n\in \mbb N}\subset X^i_x$ and $u,f\in X^i_x$. Then  $(u_n)_{n\in \mbb N}$ also converges in $W^1_1(\Omega)$ (respectively, in $C^1(\ov{\Omega})$) and hence 
$$
\lim\limits_{n\to \infty} d\Delta u_n = \lim\limits_{n\to \infty}(\mc T_0 u_n - \nabla d\cdot \nabla u_n) = f - \nabla d\cdot \nabla u
$$
 in $X^i_x$ (by the completeness of $W^1_1(\Omega)$, respectively, $C^1(\ov{\Omega})$) and $\nabla d$ being a bounded continuous function. Then, by \eqref{dass1}, $(d\Delta u_n)_{n\in \mbb N}$ converges if and only if $(\Delta u_n)_{n\in \mbb N}$ converges. Conversely, if $(\Delta u_n)_{n\in \mbb N}$ and $( u_n)_{n\in \mbb N}$ converge in $X_x$, then, similarly,  $(u_n)_{n\in \mbb N}$ converges in, respectively,  $W^1_1(\Omega)$ or $C^1(\ov{\Omega}),$ and we obtain the convergence of $(\mc T^i_0u_n)_{n\in \mbb N}$. Thus, for any $d$ satisfying \eqref{dass1}, $D(T^i_0)=D(\Delta^i)$.

 Returning to the dependence on $m$,  we denote by $T^i_{0,m}$ the realisation of the expression
$\nabla_x\cdot(d(x,m)\nabla_x\cdot)$ that generates a
semigroup in respective $X^i_x$ so that $D(T^i_{0,m})=
D(\Delta^i)$ for almost any $m$. 
 We can prove the following theorem.
\begin{theorem} Let $i\in \{0,1\}$ and 
\begin{equation}
\mbb R_+\ni m \to d(\cdot,m)\in C^1(\ov{\Omega})\;\text{be\;   Bochner-measurable}.
\label{assd}
\end{equation}
\begin{description}
\item (a) For  almost any $m\in \mbb R_+$, the operator $(T^i_{0,m},D(\Delta^i))$
generates an analytic stochastic semigroup $\sem{G_{T^i_{0,m}}}$ in $X^i_x$. 
\item (b)
The operator $\mb T^i_0$, defined in  (\ref{difeqd}), with the domain
\begin{equation}
D(\mb T^i_0) = \{\mb u\in {\mc X^i_r}:\; \mb u(\cdot,m) \in D(\Delta^i), (x,m) \to
[\mc T_0\mb u(\cdot,m)](x) \in {\mc X^i_r}\}
\label{domdif}
\end{equation}
generates an analytic substochastic semigroup  $\sem{G_{\mb T^i_0}}$ on ${\mc X^i_r}$.
\item (c) Assume that $a$ satisfies the relevant version of (\ref{Assumptionofa}). Then the operator
$(\mb T^i, D(\mb T^i)) = (\mb T^i_0+\mb A^i, D(\mb T^i_0)\cap D(\mb A^i)),$ where $\mb A^i$ was defined in \eqref{dom}, is the generator of an analytic substochastic semigroup, say $\sem{\mb{ G}_{\mb T^i}},$ and satisfies \eqref{AuEst2}. 
\end{description}
\label{bath}
\end{theorem}
\begin{proof}
Part (a) (apart from the stochasticity, which is proved
below) follows from the results discussed in the first part of the section.

To prove (b), we use Proposition \ref{propmes}, and Corollaries \ref{cormes} and \ref{ana1}. We note that the fact that we assumed in Proposition \ref{propmes} that $\Theta$ is open, while here we work in $\mbb R_+ =[0,\infty)$, does not cause any problem, as, for the measure $dm_r$, adding a single point to the domain does not change the measurability of a function.    Let $\mb u(m) = R(\la, T^i_{0,m})\mb f(m), \mb f\in \mc X^i_r$. By the first part of the proof of Proposition \ref{propmes}, $\mb u$ is measurable if $m\mapsto T^i_{0,m}u$ is continuous for any $u\in D(\Delta^i),$ for which it suffices that $d\in C^1(\mbb R_+, C^1(\ov \Omega)).$ 

Further, if $T^i_{0,m,k}$ is given by $[\mc T_{0,k} u](x,m)= \nabla_x\cdot (d_k(x,m)\nabla_x u)$ on $D(\Delta^i)$, then 
$$
\lim\limits_{k\to \infty}T^i_{0,m,k}u = T^i_{0,m}u
$$
means 
$$
\lim\limits_{k\to\infty}\left(\int\limits_\Omega |(d_k(x,m)-d(x,m))\Delta u(x)|dx + \int\limits_\Omega |\nabla_x(d_k(x,m)-d(x,m))\cdot \nabla_x u(x)|dx\right) = 0,
$$
for which it suffices that for almost every $m$, $d_k(\cdot,m)\to d(\cdot,m)$ in $C^1(\ov \Omega)$.  
Let us consider $d$ given by a step function $$d(x,m) = \sum\limits_{i=1}^n \chi_{I_i}(m) \phi_i(x),$$ for some $n \in \mbb N$, where $\phi_i\in C^1(\ov \Omega)$ and $\chi_{I_i}$ are characteristic functions of non-overlapping intervals $(a_i,b_i)\subset \mbb R_+$, $i=1,\ldots,n$. Since $d\notin C^1(\mbb R_+,C^1(\ov{\Omega}))$, we construct an approximation as follows.  By \cite[Lemma, p.8]{Vlad}, we have, for each $\chi_i,$ a sequence of $C^\infty$ functions $(\eta_{i,k})_{k\in \mbb N}$ such that $0\leq \eta_{i,k}(m)\leq 1$, $\eta_{i,k}(m)=0$ for $m \notin \left(a_i-\frac{3}{k},   b_i+\frac{3}{k}\right)$ and $\eta_{i,k}(m)=1$ for $m \in \left(a_i-\frac{1}{k},  b_i+\frac{1}{k}\right)$. Then 
$$
d_k(\cdot,m) :=  \sum\limits_{i=1}^n \eta_{i,k}(m) \phi_i(\cdot) \to \sum\limits_{i=1}^n \chi_{I_i}(m) \phi_i(\cdot) = d(\cdot,m)\;\text{as}\; k\to \infty
$$
in $C^1(\ov \Omega)$ for $m \in \mbb R_+$. Thus, $m \mapsto R(\la, T^i_{0,m})\mb f(m)$ is measurable for any $T^i_{0,m}$ if $d$ is a step $C^1(\ov\Omega)$-valued function of $m$. Finally,  by \cite[Section 1.1]{ABHN}, any $C^1(\ov\Omega)$-valued Bochner-measurable function is the limit of a sequence of step functions converging for almost any $m$ in  $C^1(\ov\Omega)$, and hence the first part of (b) follows from Corollary \ref{cormes} and Proposition \ref{glueprop1}. 

To prove that $\sem{G_{\mb T^i_0}}$ is substochastic, we note that the semigroups $\sem{G_{T^i_{0,m}}}$ with fixed $m$ are substochastic by \cite[Theorem 4.8.3]{Fa}, that is,  they are positive for  any $m$, and 
$$
\|G_{T^i_{0,m}}(t)\mb f(m)\|_{X_x^i}\leq \|\mb f(m)\|_{X_x^i}, \quad \mb f\in \mc X^i_r.
$$
Then, the integration with respect to $m$ yields the thesis for
$\sem{G_{\mb T^i_0}}$. Finally, the analyticity follows from part 1 and Corollary \ref{ana1}.

Part (c) follows directly from Proposition \ref{prop31},  Lemma \ref{lemAR}, and Corollary \ref{ana2}. 
\end{proof}
\begin{remark}
We emphasise that we do not assume that $d_{\min}$ is strictly separated from zero or that $d_{\max}$ is bounded. This is of particular importance in applications to fragmentation theory where it is expected that small particles diffuse faster than the large ones, so that we expect that $d_{\min}(m)\nearrow +\infty$ as $m\searrow  0$ and $d_{\max}(m)\searrow 0$ as $m\nearrow\infty$. 
\end{remark}
 
\subsection{Summary of the results for transport--fragmentation equation}\label{sec43}

Let us summarise the results for the full transport--fragmentation equation. 
\begin{theorem}\label{Thsum} 
\begin{description}
\item 1. Let $i\in \{0,1\},$  $\mc T_0$ be given by \eqref{ops1} with $\mb \omega$ satisfying the relevant version of (a1) and (a2), or by \eqref{difeqd} with $d$ satisfying \eqref{assd}, and  let 
$\mb T^i_0$ be the restriction of $\mc T_0$ to $D(\mb T^i_0)$ defined in \eqref{dom} for the respective $\mc T_0$. Further, let $a$ satisfy the relevant part of \eqref{Assumptionofa}.  Then, the  operator $(\mb T^i_0, D(\mb T^i_0))$ generates a substochastic semigroup (thus assumption (A1) is satisfied), and the operator $(\mb T^i, D(\mb T^i))=(\mb T^i_0+\mb A^i, D(\mb T^i_0)\cap D(\mb A^i))$ generates a substochastic semigroup with resolvent satisfying \eqref{AuEst2} in $\mc X^i_r$, for any $r\geq 0$. Furthermore, $(\mb T^i, D(\mb T^i))$ is analytic in the diffusion case. 
\item 2. In addition, assume that $b$ satisfies \eqref{b1}, \eqref{n(s)bounded}, \eqref{c_r(x)Bound}, and $r$ satisfies \eqref{r>max1l}. Then $(\mb K^1,D(\mb T^1)):=(\mb T^1+\mb B^1, D(\mb T^1))= (\mb T^1_0+\mb A^1+\mb B^1,D(\mb T^1))$ generates a positive $C_0$-semigroup, say $( G_{\mb K^1}(t))_{t\geq 0}$, on $\mathcal{X}^1_r$, which is analytic in the diffusion case. 
\item 3. In addition to assumptions of point 1., let there exist $\beta$ satisfying \eqref{c0}, and let the equi-integrability condition \eqref{betaass} be satisfied. Let $\mf K^i = \mb T^i_0+\mf A^i +\mf B^i$ be the operator defined as the restriction of the expression $\mc T_0+\mc A_1+\mc B_1$ in \eqref{TransProblem2} to $D(\mf K^i)=D(\mb T^i_0)\cap D(\mb A^i)$ and, in the case of advection, let $\mc T_0$ be independent of $m$ if $i=0$. Then there is $r_1$ such that for any $r\geq r_1$, $(\mf K^i,D(\mb T^i))$ generates a positive $C_0$-semigroup, say $(G_\mf K^i(t))_{t\geq0}$, on $\mathcal{X}^i_r$. If for $b$ in \eqref{ops2} we can find $\beta$ satisfying \eqref{bbeta} and the assumptions of this point, then also $(\mb K^i,D(\mb T^i))$ generates a positive $C_0$-semigroup on $\mc X^i_r,$  
    \begin{equation}
    G_{\mb K^i}(t)\mb f\leq  G_{\mf K^i}(t)\mb f, \quad \mb f\in \mc X^i_{r,+},
    \label{GKGK0}
    \end{equation}
    and
    \begin{equation}
    \|G_{\mb K^i}(t)\mb f\|_{\mc X^i_r}\leq  \|G_{\mf K^i}(t)\mb f\|_{\mc X^i_r} \leq M_re^{\omega_r t}\|\mb f\|_{\mc X^i_r}, \quad \mb f\in \mc X^i_r,
    \label{GKGK}
    \end{equation}
    for some constants $M_r\geq 1, \omega_r\in \mbb R.$ Furthermore, both $(G_{\mb K^i}(t))_{t\geq 0}$ and   $(G_{\mf K^i}(t))_{t\geq 0}$ are analytic in the diffusion case.
    \item 4. In the setting of point 3., the fragmentation operator $(\mf F^i, D(\mb A^i)):= (\mf A^i+\mf B^i, D(\mb A^i)), i=0,1,$ generates an analytic semigroup in $\mc X^i$. If  $\mb \omega$ (resp. $d$) is independent of $m$, then $$G_{\mf K^i}(t)\mb f = G_{\mb T^i_0}(t)(G_{\mf F^i}(t)\mb f)=G_{\mf F}(t)(G_{\mb T^i_0}(t)\mb f), \quad t\geq 0, \mb f\in \mc X^i_r,$$ and, for any $q\geq 0,$ there are  $M^{(q)}_r$ and $\omega_r$ such that         \begin{equation}\label{gfest2}
\|G_{\mb K^i}(t)\mb f\|_{\mc X^{i,(q)}_r}\leq  \|G_{\mf K^i}(t)\mb f\|_{\mc X^{i,(q)}_r} \leq  \frac{M_r^{(q)}e^{\omega_r t}}{t^q}\|\mb f\|_{\mc X^i_r},
\end{equation}
where $\omega_r$ is a constant depending on $r$ but not on $q$, and $
 \mc X_{r}^{i,(q)}$ is defined in \eqref{frps1}.
\item 5. Let $i=1$ and assume that \eqref{aCondition} is satisfied. If the assumptions of points 1. and 2. (or of point 3. that imply the former) are satisfied, then, for any $n,r$ and $p$ satisfying $\max\{1,l\}<n<p<r,$ there are constants $C>0, \theta>0$ such that 
	\begin{equation}\label{GpGr}
	|| G_{\mb K^1}(t)\mb {\mr u}||_{\mc X^1_r}\leq \frac{Ce^{\theta t}}{t^{\frac{r-n}{\gamma}}}||\mb{\mr u}||_{\mc X^1_p}, \hspace{0.5cm} \text{ for all } \mb {\mr u}\in\mathcal{X}^1_p.
	\end{equation}
\item 6. If, in the setting of points  4., \eqref{aCondition} is satisfied,  then for any $q:=r-p<\gamma,\, p\geq r_1$,
\begin{equation}\label{gfest3}
\|G_{\mb K^i}(t)\mb f\|_{\mc X^i_{r}}\leq  \| G_{\mf K^i}(t)\mb f\|_{\mc X^i_{r}} \leq  \frac{M_p^{(q)}e^{\omega_pt}}{t^\frac{q}{\gamma}}\|\mb f\|_{\mc X^i_p}.
\end{equation}
\end{description}
\end{theorem}
We note the following generalisation of \cite[Corollary 2.1]{BanLam20}. The proof in \textit{op.cit.} can be adapted to any semigroup with a moment regularising property, and thus is omitted.  
\begin{proposition}
    \label{DomainsCor}
	Let $i\in \{0,1\}$ and $\mb Q^i$ be the generator of a semigroup $\sem{G_{\mb Q^i}}$  satisfying \eqref{GTmest}, or \eqref{GpGr} (with $i=1$), or \eqref{gfest3}, in the respective setting.   Then, $G_{\mb Q^i}(t): D_p(\mb Q^i)\to D_r(\mb Q^i)$ for all $t > 0$ and hence the corresponding Cauchy problem has a classical solution in respective $\mathcal{X}^i_r$ for any $\mr{\mb u}\in \mathcal{X}^i_r\cap D_p(\mb Q^i)$.
\end{proposition} 

\subsection{Transport--fragmentation--coagulation equation}

As commented on in \cite[Section 11.2]{BLL2}, the key to dealing with fragmentation--coagulation models with spatial diffusion is ensuring that the evolution remains bounded in space. Throughout this section, we assume that $r\geq r_1,$ and $a$ satisfies \eqref{aa0c}. Then,  the transport--fragmentation semigroup behaves well in $\mc X_r^0$ for sufficiently large $r$ and has a moment regularising property there, provided the transport operator is independent of $m$ and the results of  Theorem \ref{Thsum}, points 1. and 3.--6. are valid.  
Thus, in this section, we fix $i=0$ and drop it from the notation, that is, $\mc X_r^0 =: \mc X_r$ (observe that $\mc X^0_{r}\subset  \mc X^1_r$ if $\Omega$ is bounded).  

We observe that if \eqref{aCondition} is satisfied, then, for any $0\leq q \leq \gamma,$  
$
\frac{1+m^q}{\alpha_1(m)} $ is bounded on $\mbb R_+,$  
thus $\mb A$  and $\mb A+\mb A_q$, where $\mb A_q$ is the restriction of $\mc A_q u(x,m):= -a_q(1+m^q)u(x,m)$ to $D(\mb A)$ and all results proven for operators related to $\mb T$ remain in place for $\mb T_q= \mb T+\mb A_q$. 
\subsubsection{Properties of the coagulation operator}
Assume that there are $0\leq q<\gamma$ and  $k_0$ such that for a.a. $x\in \Omega, m,s\in \mbb R_+,$ we have
\begin{equation}
0\leq k(x,m,s) \leq k_0(1+m^q)(1+s^q).
\label{kest}
\end{equation}
That the solutions to \eqref{GeneralEqn} are nonnegative is not obvious due to the presence of the negative term in the coagulation operator. Hence, we consider in $\mc X_{r}$ the modified problem 
\begin{equation}\label{rem2}
\begin{split}
\p_t u(t,x,m) &= \mc T_0u(t,x,m) +\mc Au(t,x,m) + \mc A_q u(t,x,m) +\mc Bu(t,x,m)\\
&\phantom{xx}-\mc A_q u(t,x,m) + \mc Cu(t,x,m)\\
&=: [\mb T_q\mb u](t,x,m) +[\mb B\mb u](t,x,m)+ [\mb C_{q}\mb u](t,x,m),
\end{split}
\end{equation}
	where $a_q>0$ is to be determined. As before, $\mb K_q = \mb T_q +\mb B $ generates a positive semigroup $\sem{G_{\mb K_q}}$ on $\mc X_{r}$. Moreover, by \cite[Lemma 3.9]{kerr2025}, \eqref{rem2} is equivalent to \eqref{GeneralEqn}.

The following inequalities will often be used. For  $m\geq 0,$
\begin{equation}
(1+m^\delta)\leq 2(1+m^\eta), \quad 0\leq \delta \leq \eta,\quad\textrm{and}\quad(1+m^\delta)(1+m^\eta)\leq 4(1+m^{\delta+\eta}),\quad 0\leq \delta \leq \eta.
\label{in2}
\end{equation}
Let us consider the bilinear form $\mc C_{q}$, defined by
\begin{equation}\label{mcK}
\begin{split}
[\mc C_{q}(f, g)](x,m) &:= -\mc A_q f(x,m) + \ov{\mc C}(f,g)(x,m).
\end{split}
\end{equation}
For a given $r\geq r_1$, define $r=p+q$.
\begin{proposition} \label{propCq} For any fixed $b>0$, define 
\begin{equation}
\mc U_b := \{ f \in \mc X_{r,+}:\;\| f\|_{\mc X_{r}}\leq b\}
\label{ball}
\end{equation}
and let $a_q := 2k_0b$. The operator $\mb C_{q}:\mc X_{r} \to \mc X_{p}$ is positive, bounded, globally Lipschitz continuous on $\mc U_b$ and continuously Fr\'ech{e}t differentiable as a function from $\mc X_r$ to $\mc X_p.$ 
\end{proposition}
\begin{proof}
As in \cite[Eqns (3.8)--(3.10)]{BanLam20},  for  $\mb f,\mb g \in \mc X_{r}$ and some constants $c_1,c_{2},$ we have 
\begin{equation}\label{boundedb3}
\begin{split}
\|[\mc C_{q}(\mb f,\mb g)]\|_{\mc X_{p}} &\leq  c_1 \|\mb f\|_{\mc X_{r}} + c_2\|\mb f\|_{\mc X_{r}}\|\mb g\|_{\mc X_{r}}
\end{split}
\end{equation}
 for all  $\mb f,\mb g \in \mc X_{r}.$ Hence, for $ \mb f\in \mc U_b$, by \eqref{kest}, \eqref{in2} and the fact that $q =r-p\leq   r$, 
$$
\sup\limits_{x\in \Omega}\int_{0}^{\infty} k(x,m,s)f(x,s)d s \leq 2k_0 (1 + m^q)\|\mb f\|_{\mc X_r}\leq 2k_0b(1+m^q), \quad m > 0,
$$
which leads to
\begin{equation}\label{Cpos}
(C_{q}\mb f)(x,m) \geq \frac{1}{2}\cl{0}{m}k(x,m-s,s) f(x,m-s)f(x,s)d s \geq 0.
\end{equation}
Next, using \eqref{boundedb3}, for  $\mb f \in \mc U_b$,
\begin{equation}
\|C_{q} \mb f\|_{\mc X_{p}} \leq   c_1 b + c_2b^2 =: K(\mc U),
\label{cest2a}
\end{equation}
and, using the algebraic property of bilinear forms,
$
\mc Q(\mb f,\mb f)-
\mc Q(\mb g,\mb g) = \mc Q(\mb f,\mb f-\mb g) +\mc Q(\mb f-\mb q,\mb g),
$
together with \eqref{boundedb3},  for all  $\mb f,\mb g \in \mathcal{U}_b$, we get 
\begin{equation}\label{K1B}
\begin{split}
\Vert C_{q} \mb f - C_{q} \mb g\Vert_{\mc X_{p}}  &\leq  4a_q c_1\Vert \mb f-\mb g\Vert_{\mc X_{r}} 
+ c_2 \left( \Vert \mb f \Vert_{\mc X_{r}}+ \Vert \mb g \Vert_{\mc X_{r}}\right)\Vert \mb f - \mb g \Vert_{\mc X_{r}}\\
&\leq (4a_qc_1 + 2c_2b)\Vert \mb f - \mb g \Vert_{\mc X_{r}}\leq  L(\mathcal{U}_b)\Vert \mb f - \mb g \Vert_{\mc X_{r}}. 
\end{split}
\end{equation}
The statement about the continuous Fr\'ech{e}t differentiability follows immediately from the quadratic structure of $\mb C_q$ and \eqref{boundedb3}. 
\end{proof}
\subsubsection{Solvability of the transport--fragmentation--coagulation equation}
We observe that the estimates of Proposition \ref{mcK} and \eqref{gfest3} are the only estimates used to prove local and global solvability of the fragmentation--coagulation equation with growth, \cite[Theorems 3.1 \& 3.2]{BanLam20}. Thus, a similar result can be derived here. 
\begin{theorem}\label{MT}
Let the assumptions of Theorem \ref{Thsum}, 1., 3., 4., and 6. with $i=0$ be satisfied, $p\geq r_1,$ and let \eqref{kest} hold. For any $\mb {\mr u} \in \mc X_{r,+}$ there is a mild solution to \eqref{GeneralEqn} in $\mc X_{r,+}$ defined on a maximal time interval $I_{\max}:=[0,\tau_{\mb{\mr u}}),$ and if  $\tau_{\mb{\mr u}}<0$, then $\limsup_{t\to \tau_{\mb{\mr u}}}\|\mb u(t)\|_{\mc X_r} =\infty.$

For any $\mr {\mb u}\in \mc X_r\cap D_p(\mb T)$,  see \eqref{Dr}, the mild solution is in $C(I_{\max}, \mc X_r)\cap C^1(\mr I_{\max},\mc X_r)\cap C((0,\tau_{\mb{\mr u}}),D_p(\mb T)),$ where $\mr I_{\max}=(0,\tau_{\mb{\mr u}})$, and is a classical solution to \eqref{rem2} in $\mc X_p.$
\end{theorem}
\begin{proof}
As noted above, the proof of the existence of the mild solution follows the lines of \cite[Theorem 3.1]{BanLam20}, so we only provide the opening estimates with slightly simpler proofs. 

Let $\mb {\mr u} \in \mc X_{r,+}$ be such that
\begin{equation}\label{wlfin}
\|\mb {\mr u}\|_{\mc X_{r}} \leq \frac{ b}{2}.
\end{equation}
As in, e.g., \cite[Theorem 8.1.1]{BLL2}, the mild formulation of \eqref{rem2} in $\mc X_r$ is the fixed point  problem
 \begin{equation}\label{eq3.2}
	\mb u(t) =  [\mbb F\mb u](t):=G_{\mb K_q}(t) \mr {\mb u} +  \int\limits_{0}^{t} G_{\mb K_{q}} (t - s) \mb C_{q} \mb u(s) d s,
\end{equation}
  in the space $Y_r=C([0,\tau], \mc U_b),$ with $\mc U_b$ defined by (\ref{ball}) and $\tau$ to be determined so that $\mbb F$ is a contraction on $Y_r$, when $Y_r$ is  equipped with the metric induced by the norm from $C([0,\tau], \mc X_{r})$. First, we prove that $\mbb F$ is continuous on $Y_r$.  Since $\sem{G_{\mb K_q}}$ is a $C_0$-semigroup on $\mc X_r$, we can focus on the integral term. For any $t\geq 0, 0<h <\tau-t$
 \begin{align*}
 & \left\|\int_{0}^{t+h}G_{\mb K_q}(t+h-s)\mb C_{q} \mb u(s)d s-\int_{0}^{t}G_{\mb K_q}(t-s)\mb C_{q} \mb u(s)d s\right\|_{\mc X_{r}}\\
 &\leq  \int_{0}^{h}\left\|G_{\mb K_q}(t+h-s)\mb C_{q} \mb u(s)\right\|_{\mc X_{r}}d s  
  +\int\limits_{0}^{t}\left\|G_{\mb K_q}(t-s)(\mb C_{q}\mb u(s+h)- \mb C_q\mb u(s))\right\|_{\mc X_{r}}d s \\
 &=: I_1(h)+I_2(h).
 \end{align*}
 We note that the change of variables in $I_2$ is justified as for $s\in [0,t],$ $s +h <\tau$, and $\mb u(s +h)\in \mc U_b$ for $s \in [0,t]$.
 To estimate $I_1(h)$, we first observe that for any function $\mb u\in Y_r$ and for any $\sigma>0$, $0\leq s \leq \tau$, and $q' :=\frac{q}{\gamma}$ we have, by \eqref{gfest3} and \eqref{boundedb3},
 \begin{align*}
 \left\|G_{\mb K_q}(\sigma)\mb C_{q} \mb u(s)\right\|_{\mc X_{r}} &\leq \frac{M_pe^{\omega_p\sigma}}{\sigma^{q'}}
 \left\|\mb C_{q} \mb u(s)\right\|_{\mc X_{p}}\leq \frac{M_pe^{\omega_p\sigma}}{\sigma^{q'}} K(\mc U_b).
 \end{align*}
 Thus, 
 \begin{subequations}\label{I1I2}
 \begin{equation}
 \begin{split}
 I_1(h) &\leq 
 \frac{M_pe^{\omega_p(t+h)}((t+h)^{1-q'}-t^{1-q'})}{1-q'} K(\mc U_b)\to 0\quad\text{as}\;h\to 0^+, 
 \end{split}
 \end{equation}
 and, analogously, for $t>0, t+h<\tau$, by \eqref{eq3.2},
 \begin{equation}
 \begin{split}
 I_2(h)\leq \frac{M_pe^{\omega_p t}t^{1-q'}}{1-q'} L(\mc U_b) \sup\limits_{0\leq s\leq t}\|\mb u(s+h) -\mb u(s) \|_{\mc X_{r}}.
 \end{split}
 \end{equation}
 \end{subequations}
 which converges to 0 as $h\to 0^+$, since a continuous function on a compact interval is uniformly continuous. Hence, $t\mapsto \mbb F\mb u(t)$ is right-hand continuous for $t\to \mb u(t)\in Y_r$. To prove left continuity, for any $0<t\leq \tau$ and $h>0$ such that $t-h>0$, we have 
 \begin{align*}
  &\left\|\int_{0}^{t-h}G_{\mb K_q}(t-h-s)\mb C_{q} \mb u(s)d s-\int_{0}^{t}G_{\mb K_q}(t-s)\mb C_{q} \mb u(s)d s\right\|_{\mc X_{r}}\\
 &\leq  \int_{0}^{h}\left\|G_{\mb K_q}(t-s)\mb C_{q} \mb u(s)\right\|_{\mc X_{r}}d s  
 + \int\limits_{h}^{t}\left\|G_{\mb K_q}(t-s)(\mb C_{q}(\mb u(s-h)- \mb C_q\mb u(s))\right\|_{\mc X_{r}}d s,
 \end{align*}
 which, as above, tends to zero as $h\to 0^+.$
 Now, using similar estimates, 
 \begin{align*}
\|[\mbb F]\mb u(t)\|_{\mc X_r}&\leq M_pe^{\omega_p t}\frac{b}{2} + \frac{M_pe^{\omega_p t} K(\mc U_b)}{1-q'}t^{1-q'}\\
\|[\mbb F\mb u](t)-[\mbb F\mb v](t)\|_{\mc X_r}&\leq \frac{M_pe^{\omega_p t}t^{1-q'}}{1-q'} L(\mc U_b) \sup\limits_{0\leq s\leq t}\|\mb u(s) -\mb v(s) \|_{\mc X_{r}},
 \end{align*}
   from which it easily follows that $\mbb F$ is a contraction on $Y_r$ for sufficiently small $\tau.$ The remainder of the proof is standard. 

   Similarly, the proof of the classical solvability of \eqref{rem2} is a repetition of the proof of \cite[Theorem 3.2]{BanLam20}. A significant role in the proof is played by \cite[Corollary 2.1]{BanLam20}, which was extended to the current setting in Proposition \ref{DomainsCor}. It is instructive, however, to write down the fundamental identity (correcting an editorial mistake in \textit{op.cit.}) to indicate why we only have the classical solution in a bigger space $\mc X_p$ despite $\mb u$ being differentiable in the smaller space $\mc X_r.$ To show that  a mild solution $\mb u$ satisfies $\mb u(t) \in D_p(\mb K)$ for $t>0$, we evaluate
 \begin{align*}
&	   \frac{1}{h}(G_{\mb K_q} (h) - I)\mb u(t) =
   \frac{1}{h}G_{\mb K_q} (t)(G_{\mb K_q} (h) - I) \mr{\mb u} + \frac{1}{h} G_{\mb K_q}(h)\int_{0}^{h} G_{\mb K_q} (t - s) \mb C_{q}\mb u (s)d s \\
&\phantom{xxx}- \frac{1}{h} \int_{t-h}^{t} G_{\mb K_q} (t  - s) \mb C_q\mb u(s+h)d s
	+ \frac{1}{h} \int_{0}^{t} G_{\mb K_q} (t - s)(\mb C_{q} \mb uf(s+h)-C_{q}\mb u(s)) d s \\&
\phantom{xxx} =: L_1(h) + L_2(h) + L_3(h) +L_4(h).
	\end{align*}
 Using the regularising character of $\sem{G_{\mb K_q}}$ for $t>0,$ see Proposition \ref{DomainsCor}, and the continuity of $0<t\mapsto G_{\mb K_q}(t)\mb f$ in $\mc X_r$ for any $\mb f\in \mc X_p$ (that follows by noting that $G_{\mb K_q}(t_0)\mb f\in \mc X_r$ for any $t_0>0$ and writing $G_{\mb K_q}(t)\mb f = G_{\mb K_q}(t-t_0)G_{\mb K_q}(t_0)\mb f$ for $0<t_0<t$), we establish \begin{align*}
 \lim\limits_{h\to 0^+}L_1(h) &= G_{\mb K_q}(t) \mb K_q\mr{\mb u}, \\ 
 \lim\limits_{h\to 0^+} L_2(h) &= G_{\mb K_q} (t ) \mb C_{q} \mr {\mb u}, \quad
 \lim\limits_{h\to 0^+} L_4(h) = \cl{0}{t} G_{\mb K_q}(t-s)\partial \mb C_{q}\mb u(s)\p_s\mb u(s)d s
 \end{align*}
in $\mc X_r$. However, in $L_3(h)$ we must use the continuity of the integrand at $t=0$ so we are only able to pass to the limit in $\mc X_{p}$.  Then, in the same way as for $L_2(h)$, we have
  $$
 \lim\limits_{h\to 0^+} L_3(h) =  -\mb C_{q}\mb u(t),
 $$
in $\mc X_{p}$. Hence $\mb u(t) \in D_p(\mb K)=D_p(\mb K_q)$ for $t>0$ and
 \begin{align*}
 \mb K_q\mb u(t) &= 
 - \mb C_{q} \mb u(t) + \frac{d}{dt} \mb u(t),
 \end{align*}
   where we used the integral formula for the derivative of the mild solution.  
   \end{proof}
   The restrictive assumption that the transport part $\mc T_0$ is independent of $m$ was only needed for the existence of the transport--fragmentation semigroup and the availability of  \eqref{gfest3}, necessary to prove that $G_{\mb K_q}(t)\mb C_q, t\geq 0,$ are well-defined operators on $\mc X_r$ despite $\mb C_q$ being unbounded there.   Thus, we immediately obtain
\begin{corollary}
Let $\mc T_0$ be a diffusion operator, and let the assumptions of Theorem \ref{Thsum}, points 1. and  3.  with $i=0,$ be satisfied, $r\geq r_1,$ and let  \eqref{kest} hold with $q=0,$ that is, let $k$ be bounded.  For any $\mb {\mr u} \in \mc X_{r,+}$ there is a mild solution to \eqref{GeneralEqn} in $\mc X_{r,+}$ defined on maximal time interval $[0,\tau_{\mb{\mr u}})$ and if  $\tau_{\mb{\mr u}}<0$, then $\limsup_{t\to \tau_{\mb{\mr u}}}\|\mb u(t)\|_{\mc X_r} =\infty.$

Moreover, for any $\mr {\mb u}\in  D(\mb K)\subset \mc X_r,$ the mild solution is in $C(I_{\max}, D(\mb K))\cap C^1(I_{\max},\mc X_r)$ and is a classical solution to \eqref{rem2} in $\mc X_p.$
\end{corollary}
\begin{proof}
The diffusion--fragmentation semigroup $(G_{\mb K^0}(t))_{t\geq 0}$ exists by Corollary \ref{corg1a}. Then, the proof is standard as $q= 0$ implies that $\mb A_q\mb u = k_0b\mb u$ and Proposition \ref{propCq} implies that $\mb C_q$ is globally Lipschitz in $\mc X_r$ on $\mc U_b$, and continuously Fr\'{e}chet differentiable  on $\mc X_r.$
\end{proof}

Another model to which the techniques of the proof of Theorem \ref{MT} can be used almost verbatim is  the transport--coagulation problem
\begin{equation}\label{GeneralEqn1}
\begin{split}
\partial_tu(t,x,m)&=\mathcal{T}_0u(t,x,m)+\mathcal{A}u(t,x,m)+\mathcal{C}u(t,x,m) \hspace{.5cm} t>0,~(x,m)\in\Omega\times \mathbb{R}_+,\\
u(0,x,m)&=\mr{u}(x,m), \quad (x,m)\in\Omega\times\mathbb{R}_+.
\end{split}
	\end{equation}
Here, Proposition \ref{propCq} applies unchanged and, using \eqref{GTmest} instead of \eqref{gfest3} in \eqref{I1I2}, we have
\begin{corollary}
Let $p\geq 0$ and $r=p+q$, the assumptions of Theorem \ref{Thsum}, point 1. with $i=0,$ and \eqref{aCondition}  be satisfied, and let \eqref{kest} hold. For any $\mb {\mr u} \in \mc X_{r,+}$ there is a mild solution to the problem \eqref{GeneralEqn1} in $\mc X_{r,+}$ defined on maximal time interval $[0,\tau_{\mb{\mr u}})$ and if  $\tau_{\mb{\mr u}}<0$, then $\limsup_{t\to \tau_{\mb{\mr u}}}\|\mb u(t)\|_{\mc X_r} =\infty.$

Moreover, for any $\mr {\mb u}\in \mc X_r\cap D_p(\mb T)$, the mild solution belongs to $C(I_{\max}, \mc X_r)\cap C^1(\mr I_{\max},\mc X_r)\cap C(\mr I_{\max},D_p(\mb K))$ and is a classical solution to \eqref{GeneralEqn1} in $\mc X_p.$

If $q=0$, then for any $\mr {\mb u}\in  D(\mb K)\subset \mc X_r,$ the mild solution is in $C(I_{\max}, D(\mb K))\cap C^1(I_{\max},\mc X_r)$ and is a classical solution to \eqref{GeneralEqn1} in $\mc X_p.$
\end{corollary}
\section{Conclusion}
We developed a theory of $C_0$-semigroups with a parameter and used it to provide a comprehensive theory of the spatially inhomogeneous fragmentation processes with a general transport operator in spaces $\mc X^0_r=L_1(\mbb R_+, X_x^i, (1+m^r)dm), i\in \{0,1\},$ where either $X_x^1=L_1(\Omega,dx),$ $\Omega \subseteq \mbb R^N,$ or $X_x^0$ is an appropriate space of continuous functions on $\Omega.$ Due to the unavailability of certain estimates in the latter case that makes the application of the Miyadera--Desch perturbation theorem impossible, we used a novel approach consisting in constructing a majorising $x$-independent fragmentation-like problem and showing that if the new fragmentation kernel is uniformly integrable, the transport--fragmentation problem is well-posed in $\mc X^i_r, i\in \{0,1\},$ for sufficiently large $r$ (if the advection operator is independent of $m$ in the case $i=0$), and proved that the solution semigroup has a moment regularising property if the transport operator (in both advection and diffusion case) is independent of the cluster mass $m$. We demonstrated that the theory applied to a class of advection--fragmentation and diffusion--fragmentation problems. The results for the transport--fragmentation problem are summarised in Theorem \ref{Thsum}. 

In the second part, we proved the classical local solvability of the transport--fragmentation--coagulation problem with unbounded coagulation kernels (controlled, however, by the loss term) when the transport operator was independent of $m$, or when the gain term was absent. 

Unfortunately, the standard method of proving the global solvability through deriving moment inequalities is unavailable in $\mc X^i_r$ due to the properties of the $\sup$-norm, and this question remains open; see, however, the recent paper \cite{shindin2026}, where the global solvability of the diffusion-fragmentation-coagulation equation was established in certain cases. Another area worth exploring is the analyticity of the diffusion--fragmentation semigroup, where a precise characterisation of the interpolation spaces between the domain of the generator and $\mc X_r^0$ would allow us to remove the assumption that the diffusion was independent of $m$ in the problems with fragmentation and coagulation.  The work in these directions is ongoing.\smallskip

\noindent \textbf{Conflict of interest.} The authors declare that this work does not have any conflict of interest.\smallskip

\noindent\textbf{Data availability.} Data sharing not applicable to this article as no datasets were generated or analysed during the current study.

		\bibliographystyle{plain}
\def\cprime{$'$} \def\cprime{$'$} \def\cprime{$'$} \def\cprime{$'$}
  \def\cprime{$'$} \def\cprime{$'$} \def\cprime{$'$} \def\cprime{$'$}
  \def\cprime{$'$} \def\cprime{$'$}

\end{document}